\documentclass[11pt, reqno]{amsart}

\usepackage{amsmath}
\usepackage{amscd}
\usepackage{mathrsfs}
\usepackage{amsfonts}
\usepackage{amssymb}
\usepackage{enumerate}
\usepackage{setspace}
\usepackage{dsfont}

\usepackage{array}

\usepackage{color}

\usepackage[hyperindex, colorlinks=true]{hyperref}

\usepackage[margin=1in]{geometry}

\newcommand{\eqn}{\begin{eqnarray}}
\newcommand{\een}{\end{eqnarray}}

\def\R{\Bbb R}

\def\eps{\epsilon}

\def\al{\alpha}

\newtheorem{theorem}{Theorem}[section]

\newtheorem{lemma}{Lemma}[section]
\newtheorem{prop}[theorem]{Proposition}

\theoremstyle{definition}

\newtheorem{remark}{Remark}

\numberwithin{equation}{section}

\begin{document}

\date{}
\title{Global well-posedness  for the 2D stable Muskat problem in $H^{3/2}$} 

\author{Diego C\'ordoba}
\address{Instituto de Ciencias Matem\'aticas (ICMAT),  Consejo Superior de Investigaciones Cient\'ificas, Madrid, Spain}
\email{dcg@icmat.es}

\author{Omar Lazar}
\address{Instituto de Ciencias Matem\'aticas (ICMAT),  Consejo Superior de Investigaciones Cient\'ificas, Madrid, Spain, and Departamento de An\'alisis Matem\'atico \& IMUS, Universidad de Sevilla, C/ Tarifa s/n, Campus Reina Mercedes, 41012 Sevilla, Spain}
\email{omar.lazar@icmat.es}



\maketitle
\begin{abstract}
We prove a global existence result of a unique strong solution in $\dot H^{5/2} \cap \dot H^{3/2}$ with small $\dot H^{3/2}$ semi-norm for the 2D Muskat problem, hence allowing the interface to have arbitrary large finite slopes and finite energy (thanks to the $L^{2}$ maximum principle). The proof is based on the use of a new formulation of the Muskat equation that involves oscillatory terms. Then, a careful use of interpolation inequalities in homogeneneous Besov spaces allows us to close the {\emph{a priori}} estimates.
 \end{abstract}
\bibliographystyle{plain}

\section{Introduction}

In this paper, we are interested in the Muskat problem which was introduced in \cite{Musk} by Morris Muskat in order to describe the dynamics of water and oil in sand. The Muskat problem models the motion of an interface separating two incompressible fluids in a porous medium.  One can imagine the plane $\mathbb{R}^2$ split into two regions,  say $\Gamma_1(t)$, $\Gamma_{2}(t)$ that evolve with time. We assume that the first region $\Gamma_1(t)$ is occupied by an incompressible fluid with density $\rho_1$ and the second region $\Gamma_{2}(t)$  is occupied by  another fluid with density $\rho_2$. We further assume that both fluids are immiscible. The non mixture condition allows one to consider the interface between these two fluids. This interface corresponds to their common boundary  $\partial \Gamma_1(t)$ and $\partial \Gamma_2(t)$. The velocity in each region $\Gamma_{i}(t)$ ($i=1$ or $2$) is governed by the so-called 
 Darcy's law \cite{Darcy},
  which states that the velocity depends on the gradient pressure, the gravity and the density of the fluid (which is transported by the flow) via the following relation,
\begin{equation}\label{darcy}
\frac{\mu}{\kappa}u(x,t)=-\nabla p-(0,g\rho),
\end{equation}
where $\mu$ is the constant  viscosity, $\kappa$ is the permeability of the porous media and $g$ is the gravity.  For the sake of simplicity, we may, without loss of generality, assume that all those constants are equal to 1. The system is then driven by the following transport equation 
\begin{equation}\label{transport}
\partial_{t} \rho + u\cdot\nabla \rho=0.
\end{equation}
Since the fluids are incompressible we also have
\begin{equation}\label{incom}
\nabla\cdot u=0.
\end{equation}
Equations (\ref{darcy}), (\ref{transport}) and (\ref{incom}) give rise to the so-called incompressible porous media system (IPM). Saffman and Taylor \cite{ST} pointed out that in 2D the Muskat problem is similar to the evolution of an interface in a vertical Hele-Shaw cell. 

 For the Muskat problem we can rewrite the IPM system in terms of the dynamics of the interface in between both fluids (see \cite{ambrose} and \cite{CG}). If we denote the interface by a planar curve $z(\alpha,t)$ and if  we neglect surface tension, then the interface satisfies
\begin{align*}
\partial_{t} z(\alpha,t) = \frac{\rho_2 - \rho_1}{2\pi} \int \frac{z_1(\alpha,t) - z_1(\beta,t)}{|z(\alpha,t) - z(\beta,t)|^{2}}(\partial_{\alpha}z(\alpha,t) - \partial_{\beta}z(\beta,t)) d\beta,
\end{align*}
where the curve $z$ is asymptotically flat at infinity \emph{i.e.} $(z(\al,t)-(\al,0))\rightarrow0$ as $|\al|\rightarrow\infty.$  The point $(0,\infty)$ belongs to $\Gamma_1(t)$, whereas the point $(0,-\infty)$ belongs to $\Gamma_2(t)$.  
From elementary potential theory, we can derive explicit formulas for the velocity field $u$ and the pressure $p$ from the curve $z$.
 
 A convenient way of studying the evolution of the interface is to consider this latter as a parametrized graph of a function.  When the interface is a graph of a function $z(x,t)=(x,f(x,t))$, then this characterization is preserved locally in time by the system and $f$ satisfies the contour equation
\begin{equation}
f_t(x,t)=-\rho (\Lambda f + T(f)),
\label{muskatEQ2d}
\end{equation}
where $\rho$ is equal to $\rho=\frac{\rho_2-\rho_1}{2}$, and the operator $\Lambda^{\gamma}$, $0<\gamma< 2$,  denotes the usual fractional Laplacian operator of order $\gamma$ and is defined as
$$
\Lambda^{\gamma} f =(-\Delta)^{\gamma/2} f = C_{\gamma} P.V. \int_{\mathbb R} \frac{f(x)-f(x-y)}{\vert y \vert^{1+\gamma}} dy,
$$
where $C_{\gamma}>0$ is a positive constant. In particular, when $\gamma=1$, the constant is equal to $\frac{1}{\pi}$. \\

\noindent The operator $\mathcal{H}$ denotes the Hilbert transform operator which is defined by
$$
\mathcal{H} f = \frac{1}{\pi} P.V. \int \frac{f(x-\alpha)-f(x)}{\alpha} \ d\alpha.
$$
In particular, one may easily check that $\partial_{x}\mathcal{H}=\Lambda$. \\

\noindent As for $T$, which is the nonlinear term, it is defined by
\begin{align}\label{td}
T(f)=\frac1{\pi}\int_\R\frac{\partial_{x}f(x)-\partial_{x} f(x-\alpha)}{\alpha}
\frac{\big(\frac{f(x)-f(x-\alpha)}{\alpha}\big)^2}{1+\big(\frac{f(x)-f(x-\al)}{\alpha}\big)^2}d\alpha.
\end{align}

\noindent Equivalently, the Muskat equation can be written as

\begin{equation}
\ (\mathcal{M}) \ : \\\left\{
\aligned
&\partial_{t}f=  \frac{\rho}{\pi} P.V. \ \partial_{x}\int \arctan \Delta_{\alpha} f \  d\alpha
\\ \nonumber
& f(x,0)=f_{0}(x),
\endaligned
\right.
\end{equation}
where $\Delta_{\alpha} f \equiv\frac{f(x,t)-f(x-\alpha,t)}{\alpha}=\frac{\delta_{\alpha}f(x,t)}{\alpha}.$ \\

Indeed, it is well known that linearizing $\mathcal{M}$ around the flat solution give rise the fractional heat equation $f_t=\rho \Lambda f$  (see e.g. \cite{ambrose} and \cite{CG}). The equation is linearly stable if and only if the heavier fluid is below the interface (that is $\rho_2> \rho_1$), otherwise we say that the curve is in the unstable regime. This is known as the Rayleigh-Taylor condition and is
determined by the normal component of the pressure gradient jump at the interface having a distinguished sign  (also called the Saffman-Taylor condition).

This equation has attracted the attention of the mathematical community in the past several years and we shall briefly sum up the results known regarding the Cauchy problem for  $(\mathcal{M})$ in the stable regime ($\rho >0$).  
First of all, let us recall that this equation has a maximum principle for $\|f(\cdot,t)\|_{L^\infty}$ and $\|f(\cdot,t)\|_{L^2}$.  Indeed, it is shown in \cite{CG2}  that
\begin{equation*}
\Vert f(t)||_{L^\infty(\mathbb R)}\leq \frac{\Vert f_0 \Vert_{{L^\infty}(\mathbb R)}}{1+ t}.
\end{equation*}
Moreover, the authors showed in \cite{CG2} that if $\|\partial_{x} f_0\|_{L^\infty}<1$, then $\|\partial_{x} f(\cdot,t)\|_{L^\infty}<\|\partial_{x} f_0\|_{L^\infty}$ for all $t>0$.
On the other hand, there is also an $L^{2}$ maximum principle (see \cite{CCGS}). More precisely we have
\begin{equation*}
\Vert f(T)\Vert^2_{L^2(\mathbb{R})}+\int_0^T\int_\R\int_\R \log\left(1+\left(\frac{f(\alpha,s)-f(\beta,s)}
{\alpha-\beta}\right)^2\right)d\alpha d\beta  ds =\Vert f_0\Vert^2_{L^2(\mathbb{R})}.
\end{equation*}
 This does not imply, for large initial data, a gain of regularity in the system. However, it was observed in \cite{G} that  a gain of regularity is possible even if we start with $L^2$  initial data, under the condition that the slope is initially less than 1  (see \cite{G}). Recall that the Muskat equation has a scaling: if $f$ is a solution associated to the initial data $f_0$, then the function $\lambda^{-1}f(\lambda x, \lambda t)$, $\lambda> 0$ is also a solution for the corresponding initial data $\lambda^{-1}f_{0}(\lambda x)$. In particular, the Lipschitz norm $\dot W^{1,\infty}$ is critical  as well as the homogeneous Sobolev norm $\dot H^{3/2}$. More generally, the whole family of homogeneous Besov spaces $\dot B^{1+p^{-1}}_{p,q}$ with $(p,q) \in [1,+\infty]^{2}$ is critical with respect to the scaling of $\mathcal{M}$.

As far as the local well-posedness results are concerned,  in \cite{CG}, the authors proved local existence in $H^{3}$. The authors of \cite{CGS} were able to lower the local theory to $H^{2}$.  Recently, in \cite{CGSV}, Constantin, Gancedo, Shvydkoy and Vicol have proved that the equation is locally well-posed in $\dot W^{2,p}$ with $p>1$. There is another result by Matioc \cite{Mat1} where local existence is obtained in $H^{s}$, $s \in (3/2,2]$. Instant analyticity  is obtained in \cite{CCFG1} from any initial data in $H^4$ (see also \cite{Mat1}).

If the heavier fluid is above the interface, {\emph{i.e.}}  $\rho <0$, then the equation $(\mathcal{M})$ is ill-posed in Sobolev spaces (see \cite{CG} and \cite{BSW}). However, there exists weak solutions to the (IPM) system starting with an initial data with a jump of densities  in the  unstable regime. These solutions create a zone around the initial interface where the two fluids mix. This zone grows over time, for more details see \cite{Otto}, \cite{LS}, \cite{CCF} and \cite{FL}. \\

As for the global well-posedness  results, all the available theorems were done for data having small slopes $\|\partial_{x} f_0\|_{L^\infty}<1$.  These theorems are usually obtained by taking advantage of the parabolic nature of $(\mathcal{M})$ for very small initial data. The first global well-posedness result  for small initial data is established in \cite{CG}, which was inspired by the proof in \cite{SCH} for a different setting (the fluids have same densities but different viscosities). In the critical framework, the authors of \cite{CCGS} and \cite{CCGSP} obtained a global existence result of a unique strong solution in $H^{3}$, if initially the Wiener norm (\emph{i.e.} those $f$ verifying that the Fourier transform of $\Lambda f$ is integrable) is smaller than 1/3 (which, in particular, implies a slope smaller than 1). As well, the authors of \cite{CCGS} were able to show the global propagation of the Lipschitz regularity provided the initial data has a Lipschitz norm less than 1.  Recently, the authors of \cite{CGSV}  improved the continuation criteria for solutions from the previously known $C^{2,\alpha}$ bound to instead uniformly continuous bounded slope. In \cite{Mat1}, the authors obtained the same kind of result in the space $H^{3/2+\epsilon}$ (for $\epsilon>0$) with smallness in the same space.  Furthermore, in \cite{CGSV}, the authors obtained a regularity criteria, that as long as the slope has  a uniform modulus of continuity there is existence. Recently, by taking advantage of this regularity criteria, it was proved in \cite{Cam} a global existence of smooth solutions under a smallness assumption of some critical quantity (product of sup and inf of the slope) that have to be smaller than 1. In \cite{PS} they  prove optimal decay estimates for higher derivatives in the Wiener norm. By adapting the proof of \cite{CCGS}, Deng, Lei, and  Lin \cite{DZL} have been able to prove the existence of a global solution for arbitrarily large monotonic initial data. The key ingredient in the proof is that, under a monotonicity assumption, the  slope of the interface still satisfies the $L^{\infty}$ maximum principle. Although, these solutions fail to be in $L^2$. \\

 There are also several results about the development of singularities for the Muskat equation $(\mathcal{M})$ in the stable regime. Indeed, it has been shown in \cite{CCFG1} that there exists a family of initial data (with very large slope) where the interface reaches a regime in finite time in which is no longer a graph. Therefore there exists a time $T^*$  where the slope of the solution of $(\mathcal{M})$ blows up: $\Vert f_x \Vert_{L^{\infty}}(T^*)=\infty$. In \cite{CCFG}, they proved that there exists a class of analytic initial data in the stable regime for the Muskat problem such that the solution turns to the unstable regime and no longer belongs to $C^4$. Another recent physical scenario was studied in \cite{CGZ} where the authors prove that some solutions can pass from the stable to the unstable regime and return back to the stable regime before the solution breaks down. This shift of stability phenomena illustrates the unpredictability of the solutions to the Muskat equation even starting in the stable regime. Moreover, there is numerical evidence of initial data $||\partial_{x} f_0||_{L^{\infty}}=50$ that develops an infinite slope in finite time (see \cite{CGZ1}). Similar results in a confined porous media has been obtained; for  global existence see \cite{Gr} and for shift of stability see \cite{GGr}. \\

In this paper, we develop a $\dot H^{\frac{3}{2}}$ critical theory under an arbitrary bounded slope assumption. The approach in this paper is completely new and is based on a reformulation of the usual Muskat equation $(\mathcal{M})$. This new formulation  allows one to take advantage of the oscillations which are crucial in this problem. There are many ways of measuring smoothness while trying to do {\it{a priori}}  estimates in critical spaces. Contrarily to (almost) all previous works in the Muskat equation we shall never split the study into high/low frequencies or small/big increment in the finite difference operator. On the contrary, we shall consider the interaction between both and the Besov spaces techniques will be the main tool to achieve this. It is worth saying that the new formulation of the problem turns out to be crucial to prove the main theorems of this paper since it gives new features that are very difficult to see in the original formulation $(\mathcal{M})$. If one tries to do $\dot H^{3/2}$ estimates using Besov spaces techniques by using the classical formulation of the problem as 
 stated in the introduction $(\mathcal{M})$, then one would quickly notice that there is a big issue to control the higher order terms. Therefore, a direct use of Besov space estimates would not give a satisfactory result. \\

This article is a significant step in understanding the theory of global well-posedness of large solutions in the Lipschitz space. Indeed, the main result of this paper is the global well-posedness of strong solution in $\dot H^{5/2} \cap \dot H^{3/2}$ under a smallness assumption on the $\dot H^{3/2}$ norm of the initial data. It would not be possible to prove a global results for all data in the Lipschitz space since the authors of \cite{CGZ} have shown that there is solutions with initial data having a (relatively) high slopes that become singular in finite time showing the instability of the Cauchy problem associated to initial data in critical spaces.  \\

The strategy of the proof is classical. We first establish the {\it{a priori}} estimates and then we construct a solution {\emph{via}} classical compactness arguments that allow us to pass to the weak limit in a parabolic regularized Muskat equation. \\

The outline of this paper is as follows. In the next section, we state the main results.  In the  third section, we give the definitions of the spaces along with some harmonic analysis tools that we shall use throughout the article. The fourth section is devoted to the new formulation of the problem. The fifth section is the proof of the {\it{a priori}} estimates in $\dot H^{3/2}$. In the sixth section, we give  the {\it{a priori}} estimates in $\dot H^{5/2}$. The last section is the proof of the main results.

\section{Main results}
The main result of this article is the following global existence theorem of a unique strong solution for small data in the critical Sobolev space.
\begin{theorem} \label{res}
Assume $f_0 \in \dot H^{5/2}\cap \dot H^{3/2}$  with $\Vert f_{0} \Vert_{\dot H^{3/2}}<C(\Vert \partial_{x}f_{0}\Vert_{L^{\infty}})$ small enough, then, there exists a unique strong solution $f$ to equation $(\mathcal{M})$ which verifies $ f \in L^{\infty}([0,T], \dot H^{5/2}) \cap L^{2}([0,T],  \dot H^{3})$, for all $T>0$. 
\end{theorem}
\begin{remark} For any initial data in $L^{2}$, by using the $L^{2}$ maximum principle, one gets a unique solution to equation $(\mathcal{M})$ with finite energy.
\end{remark}
\begin{remark} 
 It would be possible to lower the initial regularity by considering a data in $H^{s}$ with $3/2<s\leq 5/2$. This would lead to tedious computations, the present article does only treat the case of regular enough data with finite energy for the sake of simplicity. 
\end{remark}

In order to prove the main result, we shall prove a regularity criteria in terms of the control of the slopes that gives a condition for a weak solution to be strong.
\begin{theorem} \label{reg}
Let $f_{0} \in \dot W^{1,\infty} \cap  \dot H^{3/2}$ and let $T>0$. Assume there is a solution $f$ on $[0,T]$ with slope $\Vert \partial_{x}f \Vert_{L^{\infty}}$ and $\Vert f_{0} \Vert_{\dot H^{3/2}}<C(\Vert \partial_{x}f_{0}\Vert_{L^{\infty}})$ then there exists a unique weak solution $f$ to equation  $(\mathcal{M})$ which verifies $f \in L^{\infty}([0,T], (\dot H^{3/2}\cap \dot W^{1,\infty})) \cap L^{2}([0,T],  \dot H^{2})$ and in particular, 
 $\Vert f(T) \Vert_{ \dot H^{3/2}} \leq \Vert f_0 \Vert_{ \dot H^{3/2}}$.
\end{theorem}
\begin{remark}
 The definition of the weak solutions are easy to get. Indeed, we say that $f$ is a weak solution to the Muskat problem if, for all $\phi \in \mathcal{D}([0,T] \times \mathbb R)$, we have the following equality
 $$
 \int_{0}^{T} \int \phi_{s}(x,s) f(x,s) \ ds \ dx + \int \phi(x,0) f_{0}(x) \ dx =\int_{0}^{T} \int \phi_{x}(x,s) \int \int_{0}^{\infty}  \delta^{-1} e^{-\delta} \sin (\delta \Delta_{\alpha} f) \ d\delta \ d\alpha \ dx \ ds.
 $$
 \end{remark}

 \section{Functional setting}
 
We shall use the homogeneous Sobolev space  $\dot H^{s}$, which is endowed with the semi-norm
$$
\Vert f \Vert_{\dot H^{s}} = \Vert \Lambda^{s} f \Vert_{L^{2}}
$$

We recall the definition of the homogeneous Besov spaces  $\dot B^{s}_{p,q}(\mathbb R)$ (see \cite{B}, \cite{RS}). Let $(p,q,s) \in [1,\infty]^2 \times \mathbb R$, we say that a tempered distribution $f$ (which is such that its Fourier transform is integrable near 0)   belongs to the homogeneous Besov space $\dot B^{s}_{p,q}(\mathbb R)$ if the following quantity (which is a semi-norm), is finite, that is 
$$
\Vert f \Vert_{\dot B^{s}_{p,q}} =   \left\Vert \frac{\Vert \mathds{1}_{ ]0,1[}(s) \delta_{y}f + \mathds{1}_{[1,2[}(s) (\delta_{y}+\bar \delta_{y} f)\Vert_{L^{p}}}{ \vert{y}\vert^{s} } \right\Vert_{L^{q}(\mathbb R, \vert y \vert^{-1} dy)} <\infty,
$$
where $\delta_{y} f(x)=f(x)-f(x-y)$ and $\bar \delta_{y} f(x)=f(x)-f(x+y)$. \\

Let us recall some classical embeddings (see e.g. \cite{BL}, \cite{BCD}). We have for all $(p_1,p_2,r) \in [1,\infty]^3$
$$
\dot B^{s_1}_{p_1,r}(\mathbb R) \hookrightarrow \dot B^{s_2}_{p_2,r}(\mathbb R),
$$
where $s_1+\frac{1}{p_2}=s_2+\frac{1}{p_1}$ and $p_1<p_2$. We also have for all $(p_1, s_1) \in [2,\infty] \times \mathbb R$,
$$
\dot B^{s_1}_{p_1,r_1}(\mathbb R) \hookrightarrow \dot B^{s_1}_{p_1,r_2}(\mathbb R),
$$
for all $(r_1,r_2) \in \ ]1,\infty]^2$ such that $r_1 \leq r_2$. \\

Let $(s_1,s_2) \in \mathbb R^2$ so that $s_1<s_2$. Then, for all $\theta \in ]0,1[$ and $p \in [1,\infty]$,  we have the following real interpolation inequality
\begin{equation} \label{interp}
\Vert f \Vert_{\dot B^{\theta s_1 +(1-\theta)s_2}_{p,1}} \leq \frac{C}{s_2-s_1} \left(\frac{1}{\theta}+\frac{1}{1-\theta}\right)\Vert f \Vert^{\theta}_{\dot B^{s_1}_{p,\infty}} \Vert f \Vert^{1-\theta}_{\dot B^{s_2}_{p,\infty}}.
\end{equation}

We shall use the following generalized Calder\'on commutator type estimate which was proved by   Dawson, Macghan, and Ponce in \cite{DMP}.   Let $\Phi \in \dot W^{k+l,\infty}$ and let us consider the commutator  
$$
T_{\Phi} f= \left[\mathcal{H},\Phi \right] f,
$$
then, for all $p\in]1,\infty[$, $(k,l) \in \mathbb{N}$ and  for all $f \in L^{p}$  the following estimate holds,
\begin{equation} \label{cz}
\left\Vert T_{\Phi}  \partial^{k}_{x} f \right\Vert_{\dot W^{l,p}} \leq C_{k,l,p} \Vert \Phi \Vert_{\dot W^{k+l,\infty}} \Vert f \Vert_{L^{p}}.
\end{equation}

 Throughout the article, we shall use the notation $A \lesssim B$  if there exists a fixed constant $C>0$ such that $A \leq C B$.

\section{A new formulation of the stable Muskat equation}
 
In this section, we  give another formulation of the Muskat equation in terms of oscillatory integrals that will be useful  when doing estimates (especially to control high regularity terms) in Besov spaces.
\begin{prop} \label{NF}

Assume that $f$ solves $(\mathcal{M})$ then $f$ is a solution of $(\tilde{\mathcal{M}})$ that is
\begin{equation} 
\ (\tilde{\mathcal{M}}) \ : \\\left\{
\aligned
& f_{t} (t,x) = \frac{\rho}{\pi} P.V. \int \partial_{x}\Delta_{\alpha}f \ \int_{0}^{\infty} e^{-\delta} \cos(\delta \Delta_{\alpha} f) \ d\delta \ d\alpha \\
& f(0,x)=f_{0}(x).
\endaligned
\right.
\end{equation}

Reciprocally, if $f$ solves $(\tilde{\mathcal{M}})$ then $f$ solves $(\mathcal{M})$. That is,
$$
(\tilde{\mathcal{M}}) \iff ({\mathcal{M}}).
$$
\end{prop}

\noindent {\bf Proof  of Proposition \ref{NF}  } \\

\noindent Consider the following integrable  function
\begin{equation} 
\ \mu(x) \ = 
\aligned
\frac{1}{2}\exp(-\vert x\vert)  
\endaligned
\end{equation}

\noindent Its Fourier transform is well defined and since $\mu$ is even, we have,  

$$
\widehat{\mu} (\xi)=\int_{0}^{\infty} e^{-\delta} \cos(\delta \xi) \ d \delta.
$$
By denoting $I=\int_{0}^{\infty} e^{-\delta} \cos(\delta \xi) \ d \delta$ and by integrating by parts twice, it is not difficult to check that  $I=-\xi^2 I +1$. Therefore,

\begin{equation} \label{fourier}
\int_{0}^{\infty} e^{-\delta} \cos(\delta \xi) \ d \delta=\frac{1}{1+ \xi^{2}}.
\end{equation}

\noindent Then, by considering the restriction of this Fourier transform onto $\Delta_{\alpha} f$ one readily arrives to $(\tilde{\mathcal{M}})$. Conversely, if $f$ is a solution of $(\tilde{\mathcal{M}})$ then it is obviously also a solution of $({\mathcal{M}})$ by using once again the identity \eqref{fourier} applied to $\xi=\Delta_{\alpha}f$.

 \qed

We shall assume, without loss of generality, that $\rho=\pi$. The purpose of the next section is to prove that one has nice {\it{a priori}} estimates in $\dot H^{3/2}$.

\section{A priori estimates in $\dot H^{3/2}$} \label{h32}

We start this section by proving some useful identities that we shall use in the 
{\it{a priori}} estimates. Introduce the operator  $\bar\Delta_{\alpha}f=\frac{f(x)-f(x+\alpha)}{\alpha}$ and consider the difference $D=\Delta_{\alpha}f-\bar\Delta_{\alpha}f$. We have the following identity for $D$.
\begin{equation} \label{formuleD}
 D=\frac{f(x+\alpha)-f(x-\alpha)}{\alpha}=\frac{1}{\alpha} \int_0^\alpha f_x(x+s)+f_x(x-s)-2f_x(x) \ ds +2 f_x(x).
 \end{equation} 
We shall also need an expression for $\partial_{\alpha}D$.
\begin{eqnarray} \label{formuleDa}
\partial_{\alpha}D&=&\frac{f_x(x+\alpha)-f_x(x-\alpha)-2f_x(x)}{\alpha}-\frac{\int_{0}^{\alpha} \left(f_x(x-s)+f_x(x+s)-2f_x(x)\right) \ ds}{\alpha^2}.
\end{eqnarray}
Set $S=\Delta_{\alpha}f-\bar\Delta_{\alpha}f$. Analogously, we shall need to get some nice expression of $S$ and $\partial_{\alpha}S$. More precisely, we have the following identities.
\begin{eqnarray} \label{formuleS}
S=\Delta_{\alpha}f+\bar\Delta_{\alpha}f=-\frac{(f(x+\alpha)+f(x-\alpha)-2f(x))}{\alpha}.
\end{eqnarray} 
For $\partial_{\alpha}S$, we have that,
\begin{eqnarray} \label{formuleSa}
\partial_{\alpha}S=\bar\Delta_{\alpha} f_x-\Delta_{\alpha} f_x+\frac{f(x+\alpha)+f(x-\alpha)-2f(x)}{\alpha^2}.
\end{eqnarray}

Note that those identities will be useful since they make appear the more favorable second finite order differences. We can now state the main result of this section which is the following lemma.

\begin{lemma} \label{h3/2}
Let $T>0$, assume that $f_{0} \in \dot H^{3/2} \cap \dot W^{1,\infty}$ and that the $L^{\infty}([0,T],\dot W^{1,\infty})$ norm remains bounded. Then, if we set  $K=\displaystyle\sup_{(x,t)\in \mathbb R \times  [0,T]}\vert f_x(x,t)\vert$, we have
\begin{equation*} 
 \Vert f \Vert^{2}_{\dot H^{3/2}}(T) + \frac{\pi}{1+K^2} \int_{0}^T \Vert f \Vert^{2}_{\dot H^{2}} \ ds \lesssim \Vert f_{0} \Vert_{\dot H^{3/2}} +  \left(\Vert f \Vert_{L^{\infty}([0,T],\dot H^{3/2})} + \Vert f \Vert^{2}_{L^{\infty}([0,T],\dot H^{3/2})}  \right) \int_{0}^T \Vert f \Vert^{2}_{\dot H^{2}} \ ds.
\end{equation*}
\end{lemma}
\noindent {\bf{Proof of Lemma \ref{h3/2}}}
We multiply  $\Lambda^{3/2} f$ against $\Lambda^{3/2} f_t$  and integrate with respect to the space variable, we obtain 
\begin{eqnarray*}
\frac{1}{2}  \partial_{t} \Vert f \Vert^{2}_{\dot H^{3/2}}  &=&   \int \mathcal{H}f_{xx} \ \int \partial_{xx}  \Delta_{\alpha} f  \  \ \int_{0}^{\infty}   e^{-\delta}  \cos(\delta\Delta_{\alpha}f(x) ) \ d\delta \ d\alpha \ dx   \\
 &-&  \int  \mathcal{H}f_{xx} \ \int(\partial_{x} \Delta_{\alpha} f)^{2}  \  \ \int_{0}^{\infty} \delta  e^{-\delta}  \sin(\delta\Delta_{\alpha}f(x) ) \ d\delta \ d\alpha \ dx   \\
  &=&I_1 + I_2.
\end{eqnarray*}
Obviously, the more singular term is $I_1$. The estimation of such a term requires a long splitting into several terms that will be nicely estimated. This is the aim of the next subsection. 
\subsection*{5.1. Decomposition of the term $I_{1}$}
The aim of this subsection is to decompose $I_1$ into a sum of more easily controllable terms. The idea is to symmetrize the terms in order to make appear second finite order differences. Indeed, having a second finite order difference would allow us to share the powers of $\alpha$ in a more balanced way. Besides, estimating terms rigourously in homogeneous Besov spaces with high regularity (that is $s\in (1,2)$) requires to have second finite order differences. The main strategy will be to use the convenient identity  $\partial_{x}(f(x+\alpha)-f(x-\alpha))=\partial_{\alpha}(f(x+\alpha)+f(x-\alpha)-2f(x))$. This identity shows that $f(x+\alpha)-f(x-\alpha)=\Delta_{\alpha} f - \bar \Delta_{\alpha} f$ contains a hidden second finite order difference. To estimate the term $I_1$ we shall try to force the appearance of such symmetric terms using this idea. 

\begin{eqnarray*}
I_1&=&\int \mathcal{H}f_{xx} \ \int (\partial_{xx}  \Delta_{\alpha} f -\partial_{xx} \bar \Delta_{\alpha} f)  \  \ \int_{0}^{\infty}  e^{-\delta}  \cos(\delta\Delta_{\alpha}f(x) ) \ d\delta \ d\alpha \ dx \\
&-& \int \mathcal{H}f_{xx} \ \int \partial_{xx}  \Delta_{\alpha} f \  \ \int_{0}^{\infty}   e^{-\delta}  \cos(\delta\bar\Delta_{\alpha}f ) \ d\delta \ d\alpha \ dx \\
&=& \int \mathcal{H}f_{xx} \ \int (\partial_{xx}  \Delta_{\alpha} f -\partial_{xx} \bar \Delta_{\alpha} f)  \  \ \int_{0}^{\infty}   e^{-\delta}  \cos(\delta\Delta_{\alpha}f(x) ) \ d\delta \ d\alpha \ dx \\
&+& \int \mathcal{H}f_{xx} \ \int \partial_{xx}  \Delta_{\alpha} f \  \ \int_{0}^{\infty}   e^{-\delta} \left( \cos(\delta\Delta_{\alpha}f) -\cos(\delta\bar\Delta_{\alpha}f) \right)\ d\delta \ d\alpha \ dx \\
&-& \int \mathcal{H}f_{xx} \ \int \partial_{xx}  \Delta_{\alpha} f \  \ \int_{0}^{\infty}   e^{-\delta}  \cos(\delta\Delta_{\alpha}f) \ d\delta \ d\alpha \ dx.
\end{eqnarray*}
Hence, by symmetry (note that $\partial_{xx}\Delta_{\alpha} f -\partial_{xx} \bar \Delta_{\alpha} f$ is invariant by $\alpha \rightarrow -\alpha$), one obtains
\begin{eqnarray*}
I_1&=&\frac{1}{4} \int \mathcal{H}f_{xx}  \int (\partial_{xx}  \Delta_{\alpha} f -\partial_{xx} \bar \Delta_{\alpha} f)  \int_{0}^{\infty}  e^{-\delta} ( \cos(\delta\Delta_{\alpha}f(x))+\cos(\delta \bar\Delta_{\alpha}f(x)) \  d\delta \  d\alpha \  dx \\
&+& \frac{1}{2} \int \mathcal{H}f_{xx}  \int \partial_{xx}  \Delta_{\alpha} f \int_{0}^{\infty}   e^{-\delta} \left(  \cos(\delta\Delta_{\alpha}f) -\cos(\delta\bar\Delta_{\alpha}f)  \right)\  d\delta \ d\alpha \ dx. 
\end{eqnarray*}
By using classical trigonometry formulas, $I_1$ may be rewritten as,
\begin{eqnarray*}
I_1&=& \frac{1}{2}\int \mathcal{H}f_{xx} \int (\partial_{xx}  \Delta_{\alpha} f -\partial_{xx} \bar \Delta_{\alpha} f)  \int_{0}^{\infty}   e^{-\delta} \cos(\frac{\delta}{2}(\Delta_{\alpha}f-\bar\Delta_{\alpha}f)) \cos(\frac{\delta}{2}(\Delta_{\alpha}f+\bar\Delta_{\alpha}f)) \ d\delta \ d\alpha \  dx \\
&-& \int \mathcal{H}f_{xx} \ \int \partial_{xx}  \Delta_{\alpha} f  \int_{0}^{\infty}   e^{-\delta}  
\sin(\frac{\delta}{2}(\Delta_{\alpha}f-\bar\Delta_{\alpha}f)) \sin(\frac{\delta}{2}(\Delta_{\alpha}f+\bar\Delta_{\alpha}f))  \ d\delta \ d\alpha \ dx. 
\end{eqnarray*}
Finally, we get
\begin{eqnarray*}
I_1&=& -\int \mathcal{H}f_{xx}  \int (\partial_{xx}  \Delta_{\alpha} f -\partial_{xx} \bar \Delta_{\alpha} f)   \int_{0}^{\infty}   e^{-\delta} \cos(\frac{\delta}{2}(\Delta_{\alpha}f-\bar\Delta_{\alpha}f)) \sin^{2}(\frac{\delta}{4}(\Delta_{\alpha}f+\bar\Delta_{\alpha}f))  \ d\delta \  d\alpha \  dx \\
&+&\frac{1}{2} \int \mathcal{H}f_{xx} \ \int (\partial_{xx}  \Delta_{\alpha} f -\partial_{xx} \bar \Delta_{\alpha} f)  \  \ \int_{0}^{\infty}   e^{-\delta} \cos(\frac{\delta}{2}(\Delta_{\alpha}f-\bar\Delta_{\alpha}f))\ d\delta \ d\alpha \ dx \\
&-& \int \mathcal{H}f_{xx} \ \int \partial_{xx}  \Delta_{\alpha} f \  \ \int_{0}^{\infty}   e^{-\delta}  
\sin(\frac{\delta}{2}(\Delta_{\alpha}f-\bar\Delta_{\alpha}f)) \sin(\frac{\delta}{2}(\Delta_{\alpha}f+\bar\Delta_{\alpha}f))  \ d\delta \ d\alpha \ dx \\
&=& I_{1,1}+I_{1,2}+I_{1,3}.
\end{eqnarray*}

The aim will be to estimate the $I_{1,j}$ with $j=1,2,3$. It is important to note that the gain of regularity will come from the term $I_{1,2}$. Indeed, since $\cos(x)=1-2\sin^{2}(x/2)$, we immediatly notice that $I_{1,2}=dissipation + remainder$. The main task in this term will be to estimate the remainder. Hopefully, this remainder will be the sum of a term that will be nicely estimated plus the elliptic component (see \eqref{I121}). \\

\noindent{\bf{5.1.1. Estimates of $I_{1,1}$}} \\

\noindent To control $I_{1,1}$, we use the continuity of the Hilbert transform on $L^{2}$ along with the embedding $\dot H^{3/2} \hookrightarrow \dot B^{1}_{\infty,2}$, then one  gets
\begin{eqnarray*}
I_{1,1}&=&-\int \mathcal{H}f_{xx}  \int (\partial_{xx}  \Delta_{\alpha} f -\partial_{xx} \bar \Delta_{\alpha} f)  \int_{0}^{\infty}   e^{-\delta} \cos(\frac{\delta}{2}(\Delta_{\alpha}f-\bar\Delta_{\alpha}f)) \sin^{2}(\frac{\delta}{4}(\Delta_{\alpha}f+\bar\Delta_{\alpha}f)) \ d\delta \ d\alpha \ dx \\
&\leq&  \frac{\Gamma(3)}{4}\Vert f \Vert^{2}_{\dot H^{2}} \int \frac{\Vert \delta_{\alpha}f+\bar\delta_{\alpha}f \Vert^{2}_{L^{\infty}}}{\vert \alpha \vert^{3}} \ d\alpha \\
&\lesssim& \frac{1}{2} \Vert f \Vert^{2}_{\dot H^{2}} \Vert f \Vert^{2}_{\dot B^{1}_{\infty,2}} \\
&\lesssim& \frac{1}{2} \Vert f \Vert^{2}_{\dot H^{2}} \Vert f \Vert^{2}_{\dot H^{3/2}} \\
 \end{eqnarray*}
 
 \noindent {\bf{5.1.2. Estimates of $I_{1,2}$}}  \\
 
\noindent Recall that $I_{1,2}$ is given by 
 \begin{eqnarray*}
 I_{1,2}&=&\frac{1}{2}\int \mathcal{H}f_{xx} \ \int (\partial_{xx}  \Delta_{\alpha} f -\partial_{xx} \bar \Delta_{\alpha} f)  \  \ \int_{0}^{\infty}   e^{-\delta} \cos(\frac{\delta}{2}(\Delta_{\alpha}f-\bar\Delta_{\alpha}f))\ d\delta \ d\alpha \ dx \\
  &=& \frac{1}{2}\int \mathcal{H}f_{xx} \ \int \frac{1}{\alpha} \partial_{\alpha}\{\delta_{\alpha} f_{x} + \bar \delta_{\alpha} f_{x} \}  \  \ \int_{0}^{\infty}   e^{-\delta} \cos(\frac{\delta}{2}(\Delta_{\alpha}f-\bar\Delta_{\alpha}f))\ d\delta \ d\alpha \ dx, \\
  \end{eqnarray*}
  therefore, an integration by parts (with respect to $\alpha$) gives
   \begin{eqnarray*}
I_{1,2}  &=&\frac{1}{2}  \int \mathcal{H}f_{xx}  \int \frac{f_{x}(x-\alpha)+f_{x}(x+\alpha)-2f_{x}(x)}{\alpha^2}    \  \ \int_{0}^{\infty}   e^{-\delta} \cos(\frac{\delta}{2} (\Delta_{\alpha}f-\bar\Delta_{\alpha}f))\ d\delta \ d\alpha \ dx \\
&+&\frac{1}{4} \int \mathcal{H}f_{xx}  \int \frac{f_{x}(x-\alpha)+f_{x}(x+\alpha)-2f_{x}(x)}{\alpha}    \  \ \int_{0}^{\infty} \delta  e^{-\delta} \partial_{\alpha}\{ \Delta_{\alpha}f-\bar\Delta_{\alpha}f\}\sin(\frac{\delta}{2} (\Delta_{\alpha}f-\bar\Delta_{\alpha}f)) \\
&&\hspace{12cm}\times d\delta \ d\alpha \ dx. 
  \end{eqnarray*}
Then, in order to force the appearance of a second finite order difference, we use the identity \eqref{formuleDa}. Therefore, we get
  \begin{eqnarray*}
I_{1,2} &=& \frac{1}{2}\int \mathcal{H}f_{xx} \ \int \frac{f_{x}(x-\alpha)+f_{x}(x+\alpha)-2f_{x}(x)}{\alpha^2}    \  \ \int_{0}^{\infty}   e^{-\delta} \cos(\frac{\delta}{2} (\Delta_{\alpha}f-\bar\Delta_{\alpha}f))\ d\delta \ d\alpha \ dx \\
&+&\nonumber \frac{1}{4} \int \mathcal{H}f_{xx} \ \int \frac{f_{x}(x-\alpha)+f_{x}(x+\alpha)-2f_{x}(x)}{\alpha}    \  \ \int_{0}^{\infty} \delta  e^{-\delta} \frac{f_{x}(x-\alpha)+f_{x}(x+\alpha)-2f_{x}(x)}{\alpha} \\
&\times& \nonumber  \sin(\frac{\delta}{2} (\Delta_{\alpha}f-\bar\Delta_{\alpha}f))\ d\delta \ d\alpha \ dx \\
&-&\nonumber \frac{1}{4}\int \mathcal{H}f_{xx} \ \int \frac{f_{x}(x-\alpha)+f_{x}(x+\alpha)-2f_{x}(x)}{\alpha^3}    \  \ \int_{0}^{\infty} \delta  e^{-\delta} \int_{0}^{\alpha} f_{x}(x-s)+f_{x}(x+s)-2f_{x}(x) \ ds \\
&\times&\nonumber \sin(\frac{\delta}{2} (\Delta_{\alpha}f-\bar\Delta_{\alpha}f))\ d\delta \ d\alpha \ dx. 
\end{eqnarray*}
Then, importantly, one has to remark that the first term of $I_{1,2}$, that is \\
$$
\frac{1}{2}\int \mathcal{H}f_{xx} \ \int \frac{f_{x}(x-\alpha)+f_{x}(x+\alpha)-2f_{x}(x)}{\alpha^2}    \  \ \int_{0}^{\infty}   e^{-\delta} \cos(\frac{\delta}{2} (\Delta_{\alpha}f-\bar\Delta_{\alpha}f))\ d\delta \ d\alpha \ dx, \\
$$
contains a dissipative term. Indeed, using the trigonometry formula $\cos(\theta)=1-2\sin^{2}(\frac{1}{2}\theta)$   and the well-know identity
\begin{equation}\label{lap}
\int \frac{f_{x}(x-\alpha)+f_{x}(x+\alpha)-2f_{x}(x)}{\alpha^2} \ d\alpha=-2\pi \Lambda f_{x}=-2\pi\mathcal{H}f_{xx},
\end{equation}
one deduces that, 
\begin{eqnarray} \label{taz}
&&\frac{1}{2}\int \mathcal{H}f_{xx} \ \int \frac{f_{x}(x-\alpha)+f_{x}(x+\alpha)-2f_{x}(x)}{\alpha^2}    \int_{0}^{\infty}   e^{-\delta} \cos(\frac{\delta}{2} (\Delta_{\alpha}f-\bar\Delta_{\alpha}f))\ d\delta \ d\alpha \ dx \nonumber\\
&&= -  \pi  \int \vert \mathcal{H}f_{xx} \vert^2 \ dx - \int \mathcal{H}f_{xx}  \int \frac{f_{x}(x-\alpha)+f_{x}(x+\alpha)-2f_{x}(x)}{\alpha^2}    \int_{0}^{\infty}   e^{-\delta} \sin^{2}(\frac{\delta}{4} (\Delta_{\alpha}f-\bar\Delta_{\alpha}f))\ d\delta \ d\alpha \ dx 
\nonumber \\
&&= -  \pi \Vert f \Vert^{2}_{\dot H^{2}} - \int \mathcal{H}f_{xx}  \int \frac{f_{x}(x-\alpha)+f_{x}(x+\alpha)-2f_{x}(x)}{\alpha^2}    \int_{0}^{\infty}   e^{-\delta} \sin^{2}(\frac{\delta}{4} (\Delta_{\alpha}f-\bar\Delta_{\alpha}f))\ d\delta \ d\alpha \ dx.
\end{eqnarray}
\noindent As we shall see,  the remaining nonlinear term, namely
$$
- \int \mathcal{H}f_{xx}  \int \frac{f_{x}(x-\alpha)+f_{x}(x+\alpha)-2f_{x}(x)}{\alpha^2}    \int_{0}^{\infty}   e^{-\delta} \sin^{2}(\frac{\delta}{4} (\Delta_{\alpha}f-\bar\Delta_{\alpha}f))\ d\delta \ d\alpha \ dx,
$$
will be split into two terms. One term will be nicely estimated and the other term will be the elliptic component. Hence, using \eqref{taz}, we finally have 
\begin{eqnarray} \label{10sip}
I_{1,2}&=& \nonumber - \int \mathcal{H}f_{xx} \ \int \frac{f_{x}(x-\alpha)+f_{x}(x+\alpha)-2f_{x}(x)}{\alpha^2}    \  \ \int_{0}^{\infty}  e^{-\delta} \sin^{2}(\frac{\delta}{4} (\Delta_{\alpha}f-\bar\Delta_{\alpha}f))\ d\delta \ d\alpha \ dx \\
&+&\nonumber \frac{1}{4} \int \mathcal{H}f_{xx} \ \int \frac{f_{x}(x-\alpha)+f_{x}(x+\alpha)-2f_{x}(x)}{\alpha}    \  \ \int_{0}^{\infty} \delta  e^{-\delta} \frac{f_{x}(x-\alpha)+f_{x}(x+\alpha)-2f_{x}(x)}{\alpha} \\
&\times& \nonumber  \sin(\frac{\delta}{2} (\Delta_{\alpha}f-\bar\Delta_{\alpha}f))\ d\delta \ d\alpha \ dx \\
&-&\nonumber \frac{1}{4}\int \mathcal{H}f_{xx} \ \int \frac{f_{x}(x-\alpha)+f_{x}(x+\alpha)-2f_{x}(x)}{\alpha^3}    \  \ \int_{0}^{\infty} \delta  e^{-\delta} \int_{0}^{\alpha} f_{x}(x-s)+f_{x}(x+s)-2f_{x}(x) \ ds \\
&\times&\nonumber \sin(\frac{\delta}{2} (\Delta_{\alpha}f-\bar\Delta_{\alpha}f))\ d\delta \ d\alpha \ dx \\
&-&   \pi  \Vert f \Vert^{2}_{\dot H^{2}} \\
&=&\nonumber I_{1,2,1}+I_{1,2,2}+I_{1,2,3}+I_{1,2,4}.
 \end{eqnarray} 
 We shall estimate the $I_{1,2,j}$, $j=1,2,3$,  in the next subsection. \\
   
\noindent{\bf{5.1.2.1 Estimate of the $I_{1,2,1}$}} \\

Recall that $I_{1,2,1}$ is given by, 
 
\begin{eqnarray} \label{I121}
I_{1,2,1}&=& \nonumber - \int \mathcal{H}f_{xx} \ \int \frac{f_{x}(x-\alpha)+f_{x}(x+\alpha)-2f_{x}(x)}{\alpha^2}    \  \ \int_{0}^{\infty}   e^{-\delta} 
\sin^{2}(\frac{\delta}{4} (\Delta_{\alpha}f-\bar\Delta_{\alpha}f))\ d\delta \ d\alpha \ dx. \\
\end{eqnarray}

Since we have

\begin{equation*}
\int \frac{f_{x}(x-\alpha)+f_{x}(x+\alpha)-2f_{x}(x)}{\alpha^2} \ d\alpha=-2\pi \Lambda f_{x}=-2\pi\mathcal{H}f_{xx}.
\end{equation*}

Therefore, we notice that if we neglect the effect of the oscillatory integral we find,

\begin{equation}\label{lap}
- \int \mathcal{H}f_{xx}\int \frac{f_{x}(x-\alpha)+f_{x}(x+\alpha)-2f_{x}(x)}{\alpha^2} \ d\alpha \ dx=2\pi \int \vert \mathcal{H}f_{xx} \vert^{2} \ dx.
\end{equation}
This is a bad term since it cancels the gain of regularity obtained in the term $I_{1,2,4}$ (see \eqref{10sip}). So it is crucial to use the oscillations. The phase of the oscillatory integrals is not regular enough to be controlled in $L^{\infty}$. Indeed, estimating $\Delta_{\alpha}f-\bar\Delta_{\alpha}f$ in $L^{\infty}$ would imply a smallness on the slope. The idea is to find a second order finite difference in $\Delta_{\alpha}f-\bar\Delta_{\alpha}f$. By using the formula $\eqref{formuleD}$, we know that $\Delta_{\alpha}f-\bar\Delta_{\alpha}f$ can be written as a second finite difference plus a correction term which is $2f_{x}(x)$. The idea will be therefore to force the appearance of this correction term in the oscillatory integral. The correct quantity to substract and add in order to force the appearance of $2f_{x}(x)$  is $\sin^{2}(\frac{\delta}{2}f_{x})$. Hopefully, this quantity is independent of $\alpha$. This will allow us to remove the $\alpha$-dependence in the $\delta$-integral which is crucial in order to break the regularity issue. The remaining integral will be the desired elliptic component which will be of the form $\frac{1}{2}\frac{f^2_{x}}{1+f^{2}_{x}}$ as we shall see below. More precisely, we write

\begin{eqnarray*}
I_{1,2,1}&=&-\int \mathcal{H}f_{xx}\int \frac{f_{x}(x-\alpha)+f_{x}(x+\alpha)-2f_{x}(x)}{\alpha^2} \int_{0}^{\infty}   e^{-\delta} 
\left(\sin^{2}(\frac{\delta}{4} (\Delta_{\alpha}f-\bar\Delta_{\alpha}f))-\sin^{2}(\frac{\delta}{2}f_{x})\right)\ d\delta \ d\alpha \ dx \\
&&\ - \int \mathcal{H}f_{xx} \ \int \frac{f_{x}(x-\alpha)+f_{x}(x+\alpha)-2f_{x}(x)}{\alpha^2}    \ \int_{0}^{\infty}   e^{-\delta} 
\sin^{2}(\frac{\delta}{2}f_{x})\ d\delta \ d\alpha \ dx \\
&=&I_{1,2,1,1}+I_{1,2,1,2}.
\end{eqnarray*} 

In order to estimate $I_{1,2,1,1}$ we use the following useful identity 

\begin{eqnarray} \label{trigo}
&&\sin^{2}(\frac{\delta}{4}(\Delta_{\alpha}f-\bar \Delta_{\alpha}f)) -\sin^{2}(\frac{\delta}{2}f_x) \nonumber
= 4 \sin(\frac{\delta}{4}(\Delta_{\alpha}f-\bar \Delta_{\alpha}f-2 f_x(x))) \sin(\frac{\delta}{4}(\Delta_{\alpha}f-\bar \Delta_{\alpha}f+2 f_x(x))) \\
&&\hspace{4cm} \times  \cos(\frac{\delta}{4}(\Delta_{\alpha}f-\bar \Delta_{\alpha}f-2 f_x(x))) \cos(\frac{\delta}{4}(\Delta_{\alpha}f-\bar \Delta_{\alpha}f+2 f_x(x)))
\end{eqnarray}
So that, by using the classical estimate $\vert\sin \theta \vert \leq \vert \theta \vert$ we get 
\begin{eqnarray} \label{inegtrigo}
\left\vert  \sin^{2}(\frac{\delta}{4}(\Delta_{\alpha}f-\bar \Delta_{\alpha}f)) -\sin^{2}(\frac{\delta}{2}f_x) \right \vert \leq \delta \left\vert \Delta_{\alpha}f-\bar \Delta_{\alpha}f+2 f_x(x) \right\vert
\end{eqnarray} 
Hence, by using the identity \eqref{formuleD}  one gets
 \begin{eqnarray*}
 \vert I_{1,2,1,1} \vert 
 \lesssim \Vert f \Vert_{\dot H^{2}} \int \frac{\Vert f_{x}(x-\alpha)+f_{x}(x+\alpha)-2f_{x}(x)\Vert_{L^{\infty}}}{\vert \alpha \vert^{3}} \int_0^\alpha \frac{\Vert f_x(x+s)+f_x(x-s)-2f_x(x) \Vert_{L^{2}}}{s^q} \ s^{q}\ ds
 \end{eqnarray*}
 Let $q \in [1,\infty]$ be a real number, using H\"olderÊ\ inequality in $s$ (with $r^{-1}+ \bar r^{-1}=1$) we findÊ
 \begin{eqnarray*}
 \vert I_{1,2,1,1} \vert 
 &\lesssim& \Vert f \Vert_{\dot H^{2}} \int \frac{\Vert f_{x}(x-\alpha)+f_{x}(x+\alpha)-2f_{x}(x)\Vert_{L^{\infty}}}{\vert \alpha \vert^{3}}  \\
 &\times& \left(\int_0^\alpha \frac{\Vert f_x(x+s)+f_x(x-s)-2f_x(x) \Vert^{r}_{L^{2}}}{s^{q r}} \ ds \right)^{1/r} \left(\int_{0}^{\alpha}\ s^{q \bar r}\ ds \right)^{1/\bar r} \ d \alpha \\
 &\lesssim& \Vert f \Vert_{\dot H^{2}} \int \frac{\Vert f_{x}(x-\alpha)+f_{x}(x+\alpha)-2f_{x}(x)\Vert_{L^{\infty}}}{\vert \alpha \vert^{3-q-\frac{1}{\bar r}}}  \\
 &\times& \left(\int_0^\alpha \frac{\Vert f_x(x+s)+f_x(x-s)-2f_x(x) \Vert^{r}_{L^{2}}}{s^{q r}} \ ds \right)^{1/{r}} \ d\alpha \\
 &\leq& \Vert f \Vert_{\dot H^{2}} \Vert f_{x} \Vert_{\dot B^{2-q-\frac{1}{\bar r}}_{\infty,1}} \Vert f_{x} \Vert_{\dot B^{q-\frac{1}{r}}_{2,r}} 
  \end{eqnarray*}
 Then, by real interpolation and by choosing $q={9/8},$  $r=\bar r=2$ one obtains
  \begin{eqnarray*}
\vert I_{1,2,1,1} \vert  &\leq& \Vert f \Vert_{\dot H^{2}} \Vert f \Vert^{1/4}_{\dot B^{1}_{\infty,\infty}} \Vert f \Vert^{3/4}_{\dot B^{3/2}_{\infty,\infty}} \Vert f_{x} \Vert_{\dot B^{5/8}_{2,2}}.
  \end{eqnarray*}
 Using the interpolation inequality, we get
 $$
 \Vert f\Vert_{\dot B^{13/8}_{2,2}} \leq \Vert f\Vert^{3/4}_{\dot H^{3/2}} \Vert f\Vert^{1/4}_{\dot H^{2}},
 $$
 along with the fact that $\dot H^{1/2+\eta} \hookrightarrow \dot B^{\eta}_{\infty,\infty}$ for  $\eta=3/2$ and $\eta=1$, one finally gets
 $$
\vert I_{1,2,1,1} \vert \lesssim\Vert f \Vert^{2}_{\dot H^{2}}  \Vert f \Vert_{\dot H^{3/2}}.\\ 
$$

Then, it remains to estimate $I_{1,2,1,2}$, we write

\begin{eqnarray}
 I_{1,2,1,2} &=&\nonumber -  \int \mathcal{H}f_{xx} \ \int \frac{f_{x}(x-\alpha)+f_{x}(x+\alpha)-2f_{x}(x)}{\alpha^2}    \  \ \int_{0}^{\infty}   e^{-\delta} \ \sin^{2}( \frac{\delta}{2} f_{x}(x)) \ d\delta \ d\alpha \ dx.   \\
\end{eqnarray}
Importantly, we may explicity compute independently the $\alpha$-integral and the $\delta$-integral. Indeed, the $\alpha$-integral gives $-2\pi\mathcal{H}f_{xx}$ (see equation \eqref{lap} above). As for the $\delta$-integral, one may easily check that it is equal to
\begin{equation} \label{id}
\int_{0}^{\infty}   e^{-\delta} \ \sin^{2}( \frac{\delta}{2} f_{x}(x)) \ d\delta=\frac{1}{2}\frac{f^{2}_{x}}{1+f^{2}_{x}}.
\end{equation}
Indeed, to prove equality \eqref{id} it suffices to use equality \eqref{fourier} applied to $\xi=f_{x}(x)$.  One gets

\begin{equation*} 
\int_{0}^{\infty} e^{-\delta} \cos({\delta} f_{x}(x) ) \ d \delta=\frac{1}{1+ f_{x}(x)^{2}}.
\end{equation*}
Using the trigonometric identity $\cos(\theta)=1-2\sin^{2}(\frac{1}{2} \theta)$ applied to $\theta=\delta f_{x}(x)$ gives
\begin{equation*} 
\int_{0}^{\infty} e^{-\delta} \left(1-2\sin^{2}(\frac{\delta}{2} f_{x}(x))\right) \ d \delta=\frac{1}{1+ f_{x}(x)^{2}}.
\end{equation*}
Therefore, using that $\int_{0}^{\infty} e^{-\delta} d \delta=1$ one gets the desired identity, that is
\begin{equation*} 
\int_{0}^{\infty} e^{-\delta} \sin^{2}(\frac{\delta}{2} f_{x}(x)) \ d \delta=\frac{1}{2}\frac{f_{x}(x)^{2}}{1+ f_{x}(x)^{2}}.
\end{equation*}

\noindent Hence, combining the two identities \eqref{lap} and \eqref{id}, we get that

\begin{eqnarray} \label{a10sip}
 I_{1,2,1,5} &=& \nonumber-  \int \mathcal{H}f_{xx} \ \int \frac{f_{x}(x-\alpha)+f_{x}(x+\alpha)-2f_{x}(x)}{\alpha^2}    \  \ \int_{0}^{\infty}   e^{-\delta} \ \sin^{2}( \frac{\delta}{2} f_{x}(x)) \ d\delta \ d\alpha \ dx   \\
 &=& \pi \frac{f^{2}_{x}}{1+f^{2}_{x}} \int \vert \mathcal{H}f_{xx} \vert^{2}  \ dx.
 \end{eqnarray} 
Set  $K=\displaystyle\sup_{(x,t)\in \mathbb R \times  [0,T]}\vert f_x(x,t)\vert$. Using the fact that $x \mapsto \frac{x^2}{1+x^2}$ is an increasing function on $\mathbb R^{+}$, we get that, 
 \begin{eqnarray} \label{I1212}
 I_{1,2,1,2}  &\leq& \pi \frac{K^{2}}{1+K^{2}} \Vert f \Vert^{2}_{\dot H^{2}}.
 \end{eqnarray} 
  Therefore, we find that 

\begin{equation}
\left\vert  I_{1,2,1,1} + I_{1,2,1,2}  \right \vert \lesssim \Vert f \Vert^{2}_{\dot H^{2}} \left( \Vert f \Vert_{\dot H^{3/2}}+ \Vert f \Vert^2_{\dot H^{3/2}}\right)+   \pi\frac{K^2}{1+K^2}  \Vert f \Vert^2_{\dot H^2}
\end{equation}

\noindent {\bf{5.1.2.2. Estimate of $I_{1,2,2}$}} \\

To estimate  $I_{1,2,2}$, we observe that

\begin{eqnarray*}
I_{1,2,2} &=& \frac{1}{4}\int \mathcal{H}f_{xx} \ \int \frac{f_{x}(x-\alpha)+f_{x}(x+\alpha)-2f_{x}(x)}{\alpha}    \  \ \int_{0}^{\infty} \delta  e^{-\delta} \frac{f_{x}(x-\alpha)+f_{x}(x+\alpha)-2f_{x}(x)}{\alpha} \\
&&  \hspace{2cm} \times \sin(\frac{\delta}{2} (\Delta_{\alpha}f-\bar\Delta_{\alpha}f))\ d\delta \ d\alpha \ dx \\
&\lesssim& \Vert f \Vert_{\dot H^{2}} \left(\int 
\frac{\Vert f_{x}(x-\alpha)+f_{x}(x+\alpha)-2f_{x}(x)\Vert^{2}_{L^{2}}}{\alpha^{2}} \ d\alpha\right)^{1/2} \\
&& \hspace{3cm} \times \left(\int 
\frac{\Vert f_{x}(x-\alpha)+f_{x}(x+\alpha)-2f_{x}(x)\Vert^{2}_{L^{\infty}}}{\alpha^{2}} \ d\alpha\right)^{1/2} \\
&\lesssim& \Vert f \Vert_{\dot H^{2}} \Vert f_{x} \Vert_{\dot B^{1/2}_{\infty,2}} \Vert f \Vert_{\dot H^{3/2}} \\
&\lesssim& \Vert f \Vert^{2}_{\dot H^{2}} \Vert f \Vert_{\dot H^{3/2}}
\end{eqnarray*}
To estimate $I_{1,2,3}$ we write
\begin{eqnarray*}
I_{1,2,3}&=&-\frac{1}{4}\int \mathcal{H}f_{xx} \ \int \frac{f_{x}(x-\alpha)+f_{x}(x+\alpha)-2f_{x}(x)}{\alpha^3}    \  \ \int_{0}^{\infty} \delta  e^{-\delta} \int_{0}^{\alpha} f_{x}(x-s)+f_{x}(x+s)-2f_{x}(x)  \ ds \\
&\times& \sin(\frac{\delta}{2} (\Delta_{\alpha}f-\bar\Delta_{\alpha}f))\ d\delta \ d\alpha \ dx. \\
\end{eqnarray*}
Hence,
\begin{eqnarray*}
I_{1,2,3}&\leq&\frac{1}{4}\Vert f \Vert_{\dot H^{2}} \int\int_{0}^{\infty} \delta e^{-\delta} \frac{\Vert f_{x}(x-\alpha)+f_{x}(x+\alpha)-2f_{x}(x) \Vert_{L^{\infty}}}{\vert \alpha \vert^{3}} \\
 && \ \times \int_{0}^{\alpha}  {\Vert f_{x}(x-s)-f_{x}(x) \Vert_{L^{2}}}+{\Vert f_{x}(x+s)-f_{x}(x) \Vert_{L^{2}}} \ ds  \ d \delta \ d\alpha \\
 &\lesssim& \Vert f \Vert_{\dot H^{2}} \left( \int \frac{\Vert f_{x}(x-\alpha)+ f_{x}(x+\alpha)-2f_{x}(x) \Vert_{L^{\infty}}}{\vert \alpha \vert^{3-q-\frac{1}{\bar r}}}  \ d\alpha \right)  \\
  && \ \times \left[\left(\int  \frac{\Vert f_{x}(x-s)-f_{x}(x) \Vert^{r} _{L^{2}}}{\vert s \vert^{qr}} \ ds\right)^{1/r} + \left(\int  \frac{\Vert f_{x}(x+s)-f_{x}(x) \Vert^{r} _{L^{2}}}{\vert s \vert^{qr}} \ ds\right)^{1/r}\right]   \\
   &\lesssim&  \Vert f \Vert_{H^{2}} \Vert f_{x} \Vert_{\dot B^{2-q-\frac{1}{\bar r}}_{\infty,1}} \Vert f_{x} \Vert_{\dot B^{q-\frac{1}{r}}_{2,r}} \\
\end{eqnarray*}
Then, by choosing  $q=9/8,$ $r=\bar r=2$, and real interpolation (to get a control of $\dot B^{11/8}_{\infty,1}$) along with  classical homogeneous Besov embeddings  one gets
\begin{eqnarray*}
\vert I_{1,2,3} \vert &\leq& \frac{1}{2} \Vert f \Vert_{\dot H^{2}} \Vert f\Vert_{\dot B^{13/8}_{2,2}}  \Vert f\Vert_{\dot B^{11/8}_{\infty,1}}    \\
&\leq& \frac{1}{2} \Vert f \Vert_{\dot H^{2}}  \Vert f \Vert^{1/4}_{\dot B^{1}_{\infty,1}} \Vert f \Vert^{3/4}_{\dot B^{3/2}_{\infty,1}}
 \Vert f \Vert^{1/4}_{\dot B^{2}_{2,2}} \Vert f \Vert^{3/4}_{\dot B^{2}_{2,2}}
\\
&\leq& \frac{1}{2} \Vert f \Vert^{2}_{\dot H^{2}}  \Vert f \Vert_{\dot H^{3/2}}.
\end{eqnarray*}
Recalling that \eqref{10sip} is a dissipative term, and that \eqref{a10sip} is the elliptic component, we have proved that
 \begin{equation*} 
  \vert I_{1,2} \vert \lesssim \Vert f \Vert^{2}_{H^{2}}   (\Vert f \Vert^{2}_{\dot H^{3/2}} + \Vert f \Vert_{\dot H^{3/2}})- \pi \Vert f \Vert^{2}_{\dot H^{2}}+ \pi\frac{K^2}{1+K^2} \Vert f \Vert^{2}_{\dot H^{2}},
 \end{equation*}
where  $K=\displaystyle\sup_{(x,t)\in \mathbb R \times  [0,T]}\vert f_x(x,t)\vert$. Finally,
  \begin{equation} \label{est1}
 \vert I_{1,2} \vert \lesssim \Vert f \Vert^{2}_{H^{2}}   (\Vert f \Vert^{2}_{\dot H^{3/2}} + \Vert f \Vert_{\dot H^{3/2}})- \frac{\pi}{1+K^2} \Vert f \Vert^{2}_{\dot H^{2}}.
 \end{equation}

 \noindent{\bf{5.1.3. Estimates of $I_{1,3}$}}  \\

Let us estimate $I_{1,3}$. Let us recall that

\begin{eqnarray*}
\displaystyle I_{1,3}&=& - \int \mathcal{H}f_{xx} \ \int \partial_{xx}  \Delta_{\alpha} f \  \ \int_{0}^{\infty}   e^{-\delta}  
\sin(\frac{\delta}{2}(\Delta_{\alpha}f-\bar\Delta_{\alpha}f)) \sin(\frac{\delta}{2}(\Delta_{\alpha}f+\bar\Delta_{\alpha}f))  \ d\delta \ d\alpha \ dx.  \\
\end{eqnarray*}
The estimate of $I_{1,3}$ is not immediate. Indeed, if we do a direct estimate of $I_{1,3}$ we get that
\begin{eqnarray*}
I_{1,3}&\lesssim& \Vert f \Vert^{2}_{\dot H^{2}} \int \frac{\Vert \delta_{\alpha}f + \bar \delta_{\alpha}f \Vert_{L^{\infty}}}{\alpha^2} d \alpha = \Vert f \Vert^{2}_{\dot H^{2}} \Vert f \Vert_{\dot B^{1}_{\infty,1}}
\end{eqnarray*}
Using the interpolation inequality \eqref{interp}, for $p=\infty$ and $\theta s_1 +(1-\theta)s_2=1$ where $s_1<s_2$ one obtains
\begin{equation}  \label{interp2}
\Vert  f \Vert_{\dot B^{1}_{\infty,1}} \lesssim  \Vert f \Vert^{\theta}_{\dot B^{s_1}_{\infty,\infty}} \Vert f \Vert^{1-\theta}_{\dot B^{s_2}_{\infty,\infty}}.
\end{equation}
The conditions $\theta s_1 +(1-\theta)s_2=1$ and $s_1<s_2$ necessarily imply that $s_2>1$. Therefore, this would give a control of  $\Vert  f \Vert_{\dot B^{1}_{\infty,1}}$ by a quantity of the type $\Vert  f \Vert_{\dot B^{1+\epsilon}_{\infty,\infty}}$ which is a subcritical norm with respect to the scaling of the equation (recall that the critical space in that scale is $\dot B^{1}_{\infty,\infty})$. \\

 Since there is a regularity issue, the idea is to try to balance the derivatives. The idea is to write that $f_{xx}(x-\alpha)-f_{xx}(x)=\partial_{\alpha}(f_{x}(x)-f_{x}(x-\alpha))-f_{xx}(x)$. This kind of splitting is usually not very helpful because the finite difference is always more regular than each term written separatly. However, the presence of $\mathcal{H}f_{xx}$ in front of $f_{xx}(x)$  creates a new regularizing effect since it has a nice commutator structure. More precisely, we write

\begin{eqnarray*}
I_{1,3}&=&- \int \mathcal{H}f_{xx} \ \int \frac{f_{xx}(x)-f_{xx}(x-\alpha)}{\alpha} \  \ \int_{0}^{\infty}   e^{-\delta}  
\sin(\frac{\delta}{2}(\Delta_{\alpha}f-\bar\Delta_{\alpha}f)) \sin(\frac{\delta}{2}(\Delta_{\alpha}f+\bar\Delta_{\alpha}f))  \ d\delta \ d\alpha \ dx \\
&=& \int \mathcal{H}f_{xx} \ \int \frac{f_{xx}(x-\alpha)-f_{xx}(x)}{\alpha} \  \ \int_{0}^{\infty}   e^{-\delta}  
\sin(\frac{\delta}{2}(\Delta_{\alpha}f-\bar\Delta_{\alpha}f)) \sin(\frac{\delta}{2}(\Delta_{\alpha}f+\bar\Delta_{\alpha}f))  \ d\delta \ d\alpha \ dx \\
&=& \int \mathcal{H}f_{xx} \ \int \frac{\partial_{\alpha}\left\{f_{x}(x)-f_{x}(x-\alpha)\right\}}{\alpha}  \int_{0}^{\infty}   e^{-\delta}  
\sin(\frac{\delta}{2}(\Delta_{\alpha}f-\bar\Delta_{\alpha}f)) \sin(\frac{\delta}{2}(\Delta_{\alpha}f+\bar\Delta_{\alpha}f))  \ d\delta \ d\alpha \ dx \\
&-& \underbrace{\int \mathcal{H}f_{xx} \ \int \frac{f_{xx}(x)
}{\alpha} \int_{0}^{\infty}  e^{-\delta}  
\sin(\frac{\delta}{2}(\Delta_{\alpha}f-\bar\Delta_{\alpha}f)) \sin(\frac{\delta}{2}(\Delta_{\alpha}f+\bar\Delta_{\alpha}f))  \ d\delta \ d\alpha \ dx.}_{commutator \ structure} 
\end{eqnarray*}
By integrating by parts (in $\alpha$) and by using the notations $D=\Delta_{\alpha}f-\bar\Delta_{\alpha}f$ and $S=\Delta_{\alpha}f+\bar\Delta_{\alpha}f$, one finds
\begin{eqnarray} \label{I13}
I_{1,3}&=& \int \mathcal{H}f_{xx} \int \frac{f_{x}(x)-f_{x}(x-\alpha)}{\alpha^2} \int_{0}^{\infty}   e^{-\delta}  
\sin(\frac{\delta}{2}(\Delta_{\alpha}f-\bar\Delta_{\alpha}f)) \sin(\frac{\delta}{2}(\Delta_{\alpha}f+\bar\Delta_{\alpha}f))  \ d\delta \ d\alpha \ dx \\ \nonumber
&-& \frac{1}{2}  \int \mathcal{H}f_{xx}  \int \frac{f_{x}(x)-f_{x}(x-\alpha)}{\alpha} \int_{0}^{\infty} \delta  e^{-\delta}  
\partial_{\alpha} D \cos(\frac{\delta}{2}(\Delta_{\alpha}f-\bar\Delta_{\alpha}f)) \sin(\frac{\delta}{2}(\Delta_{\alpha}f+\bar\Delta_{\alpha}f))  \ d\delta \ d\alpha \ dx \\ \nonumber
&-&\frac{1}{2}  \int \mathcal{H}f_{xx}  \int \frac{f_{x}(x)-f_{x}(x-\alpha)}{\alpha}  \int_{0}^{\infty} \delta  e^{-\delta}  
 \partial_{\alpha} S \sin(\frac{\delta}{2}(\Delta_{\alpha}f-\bar\Delta_{\alpha}f))  \cos(\frac{\delta}{2}(\Delta_{\alpha}f+\bar\Delta_{\alpha}f))  \ d\delta \ d\alpha \ dx \\ \nonumber
 &-& \int \mathcal{H}f_{xx} \int \frac{f_{xx}(x)}{\alpha} \int_{0}^{\infty}   e^{-\delta}  
\sin(\frac{\delta}{2}(\Delta_{\alpha}f-\bar\Delta_{\alpha}f)) \sin(\frac{\delta}{2}(\Delta_{\alpha}f+\bar\Delta_{\alpha}f))  \ d\delta \ d\alpha \ dx  \\ \nonumber
&=& \sum_{i=1}^4 I_{1,3,i} 
\end{eqnarray}

 In the next subsection, we shall estimate the $I_{1,3,i}$ for $i=1,..,4$.  \\

 \noindent{\bf{5.1.3.1. Estimates of $I_{1,3,1}$}} \\
 
This term is estimated as follows,  by observing that $\dot H^{2} \hookrightarrow \dot B^{3/2}_{\infty,2}$, we may for instance write
 \begin{eqnarray*}
\vert I_{1,3,1} \vert &\leq&  \Vert f \Vert_{\dot H^{2}}  \  \ \int_{0}^{\infty} \delta  e^{-\delta} \frac{\Vert f_{x}(x)-f_{x}(x-\alpha)\Vert_{L^{2}}\Vert f(x-\alpha)+f(x+\alpha)-2f(x)\Vert_{L^{\infty}}}{\vert \alpha \vert^{3}} \ d\delta \ d\alpha  \\
&\leq&  \Gamma(2) \Vert f \Vert_{\dot H^{2}} \left(\int \frac{\Vert f_{x}(x)-f_{x}(x-\alpha) \Vert^{2}_{L^{2}}}{\vert \alpha \vert^2} \ d\alpha \int \frac{\Vert f(x-\alpha)+f(x+\alpha)-2f(x) \Vert^{2}_{L^{\infty}}}{\vert \alpha \vert^4} \ d\alpha \right)^{1/2} \\
&\lesssim& \Vert f \Vert_{\dot H^{2}} \Vert f_{x} \Vert_{\dot B^{1/2}_{2,2}} \Vert f \Vert_{\dot B^{3/2}_{\infty,2}} \\
&\lesssim&  \Vert f \Vert^2_{\dot H^{2}} \Vert f \Vert_{\dot H^{3/2}}
 \end{eqnarray*}
  
  \noindent{\bf{5.1.3.2. Estimates of $I_{1,3,2}$}} \\
  
Recall that,
\begin{eqnarray*}
I_{1,3,2}=- \frac{1}{2}  \int \mathcal{H}f_{xx}  \int \frac{f_{x}(x)-f_{x}(x-\alpha)}{\alpha} \int_{0}^{\infty} \delta  e^{-\delta}  
\partial_{\alpha} D \cos(\frac{\delta}{2}(\Delta_{\alpha}f-\bar\Delta_{\alpha}f)) \sin(\frac{\delta}{2}(\Delta_{\alpha}f+\bar\Delta_{\alpha}f))  \ d\delta \ d\alpha \ dx.
\end{eqnarray*}
Hence, by using the identity \eqref{formuleDa}, one gets
and we infer that, 
\begin{eqnarray} \label{I132}
I_{1,3,2}&=& -\frac{1}{2} \int \mathcal{H}f_{xx} \ \int \frac{f_{x}(x)-f_{x}(x-\alpha)}{\alpha} \  \ \int_{0}^{\infty} \delta  e^{-\delta}  
\frac{f_x(x+\alpha)+f_x(x-\alpha)-2f_x(x)}{\alpha} \\
&& \times \ \cos(\frac{\delta}{2}(\Delta_{\alpha}f-\bar\Delta_{\alpha}f)) \sin(\frac{\delta}{2}(\Delta_{\alpha}f+\bar\Delta_{\alpha}f))  \ d\delta \ d\alpha \ dx \nonumber \\
&+&\frac{1}{2}  \int \mathcal{H}f_{xx} \ \int \frac{f_{x}(x)-f_{x}(x-\alpha)}{\alpha} \  \ \int_{0}^{\infty} \delta  e^{-\delta}  
 \ \frac{\int_{0}^{\alpha} \left(f_x(x-s)+f_x(x+s)-2f_x(x)\right) \ ds}{\alpha^2} \nonumber \\
 &&\times \cos(\frac{\delta}{2}(\Delta_{\alpha}f-\bar\Delta_{\alpha}f)) \sin(\frac{\delta}{2}(\Delta_{\alpha}f+\bar\Delta_{\alpha}f))  \ d\delta \ d\alpha \ dx \nonumber \\ 
 &=& I_{1,3,2,1} + I_{1,3,2,2}.
\end{eqnarray}

We may estimate those two terms as follows. For the first one, it suffices to use the embedding   $\dot H^1 \hookrightarrow \dot B^{1/2}_{\infty,2}$, indeed we have
\begin{eqnarray*}
\vert  I_{1,3,2,1} \vert &\leq& \frac{1}{2}\Vert f \Vert_{\dot H^{2}} \int \frac{\Vert f_{x}(x)-f_{x}(x-\alpha) \Vert_{L^{2}}}{\vert\alpha\vert} \   
\frac{\Vert f_{x}(x+\alpha)+f_{x}(x-\alpha)-2f_{x}(x) \Vert_{L^{\infty}}}{\vert\alpha\vert} \ d\alpha  \\
&\lesssim&  \Vert f \Vert_{\dot H^{2}} \Vert f_{x} \Vert_{\dot B^{1/2}_{2,2}} \Vert f_{x} \Vert_{\dot B^{1/2}_{\infty,2}} \\
&\lesssim& \Vert f \Vert^{2}_{\dot H^{2}} \Vert f\Vert_{\dot H^{3/2}}
\end{eqnarray*}
For   $I_{1,3,2,2}$, for some $q$, $r$ and $\bar r$ (so that $1/r+1/\bar r=1$) that will be chosen latter, we write
\begin{eqnarray*}
\vert  I_{1,3,2,2} \vert &\leq&\frac{1}{2}\Vert f \Vert_{\dot H^{2}} \ \int \frac{\Vert f_{x}(x)-f_{x}(x-\alpha) \Vert_{L^{\infty}}}{\vert \alpha \vert^{3}} \  \ \int_{0}^{\infty} \delta  e^{-\delta} \\
&\times&   \vert \alpha \vert^{q+\frac{1}{\bar r}} \left(\int_{0}^{\alpha} \frac{ \Vert f_x(x-s)+f_x(x+s)-2f_x(x) \Vert^{r}_{L^{2}}}{s^{qr}} \ ds\right)^{1/r}\\
 &&\times  \frac{\Vert \delta_{\alpha}f+\bar\delta_{\alpha}f \Vert_{L^{\infty}}}{\vert\alpha \vert}  \ d\delta \ d\alpha  \\
 &\leq&\frac{\Gamma(2)}{2}\Vert f \Vert_{\dot H^{2}} \Vert f_{x} \Vert_{ \dot B^{q-\frac{1}{r}}_{2,r}} \ \int \frac{\Vert f_{x}(x)-f_{x}(x-\alpha) \Vert_{L^{\infty}}}{\vert \alpha \vert^{3-q-\frac{1}{\bar r}}} \  \frac{\Vert \delta_{\alpha}f+\bar\delta_{\alpha}f \Vert_{L^{\infty}}}{\vert\alpha \vert} \ d\alpha  \\
 &\lesssim& \Vert f \Vert_{\dot H^{2}} \Vert f\Vert_{ \dot B^{q+1-\frac{1}{r}}_{2,r}} \left(\int \frac{\Vert f_{x}(x)-f_{x}(x-\alpha) \Vert^{2}_{L^{\infty}}}{\vert \alpha \vert^{4-2q-\frac{2}{\bar r}}} \ d\alpha \int \frac{\Vert \delta_{\alpha}f+\bar\delta_{\alpha}f \Vert^{2}_{L^{\infty}}}{\alpha^{4}} \ d\alpha \right)^{1/2} \\
 &\lesssim&\Vert f \Vert_{\dot H^{2}} \Vert f \Vert_{ \dot B^{q+1-\frac{1}{r}}_{2,r}} \Vert f \Vert_{\dot B^{\frac{5}{2}-q-\frac{1}{\bar r}}_{\infty,2}}\Vert f \Vert_{\dot B^{3/2}_{\infty,2}} \\
 &\lesssim& \Vert f \Vert_{\dot H^{2}} \Vert f \Vert_{ \dot B^{q+\frac{1}{2}}_{2,r}} \Vert f \Vert_{\dot B^{2-q}_{\infty,2}}\Vert f \Vert_{\dot B^{3/2}_{\infty,2}}
 \end{eqnarray*}
Then, by choosing $\bar r= r=2$ and $q=1$, we obtain
\begin{eqnarray*}
\vert  I_{1,3,2,2} \vert &\lesssim&  \Vert f \Vert^{2}_{\dot H^{2}} \Vert f \Vert^{2}_{\dot H^{3/2}}
\end{eqnarray*}
Finally, we have shown that
\begin{eqnarray*} \vert I_{1,3,2} \vert \lesssim \Vert f \Vert^{2}_{\dot H^{2}} \Vert f \Vert^{2}_{\dot H^{3/2}}
\end{eqnarray*}

  \noindent{\bf{5.1.3.3. Estimates of $I_{1,3,3}$}} \\

\noindent Recall that
\begin{eqnarray*}
  I_{1,3,3} &=&-\frac{1}{2}\int \mathcal{H}f_{xx}  \int \frac{f_{x}(x)-f_{x}(x-\alpha)}{\alpha} \int_{0}^{\infty} \delta  e^{-\delta}   \sin(\frac{\delta}{2}(\Delta_{\alpha}f-\bar\Delta_{\alpha}f)) \partial_{\alpha} S  \cos(\frac{\delta}{2}(\Delta_{\alpha}f+\bar\Delta_{\alpha}f))  \\
  && \hspace{13cm} \times \  d\delta \ d\alpha \ dx 
  \end{eqnarray*}
  By using the identity \eqref{formuleSa}, one finds
  \begin{eqnarray*}
  I_{1,3,3}&=&  -\frac{1}{2}\int \mathcal{H}f_{xx}  \int \frac{f_{x}(x)-f_{x}(x-\alpha)}{\alpha} \ \int_{0}^{\infty} \delta  e^{-\delta}  
 \ \sin(\frac{\delta}{2}(\Delta_{\alpha}f-\bar\Delta_{\alpha}f)) \bar \Delta_{\alpha} f_{x} \ \cos(\frac{\delta}{2}(\Delta_{\alpha}f+\bar\Delta_{\alpha}f))  \\ 
 && \hspace{13cm} \times \  d\delta \ d\alpha \ dx \\
&+&\frac{1}{2} \int \mathcal{H}f_{xx}  \int \frac{f_{x}(x)-f_{x}(x-\alpha)}{\alpha}   \ \int_{0}^{\infty} \delta  e^{-\delta} 
 \ \sin(\frac{\delta}{2}(\Delta_{\alpha}f-\bar\Delta_{\alpha}f)) \Delta_{\alpha}f_{x} \ \cos(\frac{\delta}{2}(\Delta_{\alpha}f+\bar\Delta_{\alpha}f))  \\ 
 && \hspace{13cm} \times \  d\delta \ d\alpha \ dx \\ &-&\frac{1}{2} \int \mathcal{H}f_{xx} \ \int \frac{f_{x}(x)-f_{x}(x-\alpha)}{\alpha}  \int_{0}^{\infty} \delta  e^{-\delta}   \sin(\frac{\delta}{2}(\Delta_{\alpha}f-\bar\Delta_{\alpha}f)) \frac{f(x+\alpha)+f(x-\alpha)-2f(x)}{\alpha^2} \\
  && \times\cos(\frac{\delta}{2}(\Delta_{\alpha}f+\bar\Delta_{\alpha}f))  \ d\delta \ d\alpha \ dx \\
  &=& \sum_{i=1}^{3}    I_{1,3,3, i} 
  \end{eqnarray*}
    
  \noindent The control of $I_{1,3,3,1}$ is relatively easy, indeed, it suffices to write
  
  \begin{eqnarray*}
  \vert I_{1,3,3,1} \vert &\leq&\frac{\Gamma(2)}{2} \Vert f \Vert_{\dot H^{2}} \left(\int \frac{\Vert f_{x}(x)-f_{x}(x-\alpha) \Vert^{2}_{L^{4}}}{\vert \alpha \vert^{2}} \ d\alpha \int \frac{\Vert f_{x}(x)-f_{x}(x+\alpha) \Vert^{2}_{L^{4}}}{\vert \alpha \vert^{2}} \ d\alpha \right)^{1/2} \\
   &\lesssim& \Vert f \Vert_{\dot H^{2}} \Vert f_{x} \Vert^{2}_{\dot B^{1/2}_{4,2}} \\
    & \lesssim&  \Vert f \Vert_{\dot H^{2}} \Vert f \Vert^{2}_{\dot B^{7/4}_{2,2}} \\
    & \lesssim&  \Vert f \Vert^{2}_{\dot H^{2}} \Vert f \Vert_{\dot H^{3/2}}
  \end{eqnarray*}
  As well, one may easily estimate $I_{1,3,3,2}$ by writting
   \begin{eqnarray*}
   \vert I_{1,3,3,2} \vert  &\leq&\frac{\Gamma(2)}{2} \Vert f \Vert_{\dot H^{2}} \int \frac{\Vert f_{x}(x)-f_{x}(x-\alpha) \Vert^{2}_{L^{4}}}{\vert \alpha \vert^{2}} \ d\alpha  \\
   & \lesssim& \Vert f \Vert_{\dot H^{2}} \Vert f_{x} \Vert^{2}_{\dot B^{1/2}_{4,2}} \\
    & \lesssim& \Vert f \Vert_{\dot H^{2}} \Vert f \Vert^{2}_{\dot B^{7/4}_{2,2}} \\
    & \lesssim& \Vert f \Vert^{2}_{\dot H^{2}} \Vert f \Vert_{\dot H^{3/2}}
  \end{eqnarray*}
  For $I_{1,3,3,3}$, it suffices to see that
  \begin{eqnarray*}
   \vert I_{1,3,3,3} \vert &\leq&\frac{\Gamma(2)}{2} \Vert f \Vert_{\dot H^{2}} \left(\int \frac{\Vert f_{x}(x)-f_{x}(x-\alpha) \Vert^{2}_{L^{4}}}{\vert \alpha \vert^{2}} \ d\alpha \int \frac{\Vert f(x+\alpha)+f(x-\alpha)-2f(x)\Vert^{2}_{L^{4}}}{\vert \alpha \vert^{4}} \ d\alpha \right)^{1/2} \\
   & \lesssim&\Vert f \Vert_{\dot H^{2}} \left(\int \frac{\Vert f_{x}(x)-f_{x}(x-\alpha) \Vert^{2}_{L^{4}}}{\vert \alpha \vert^{2}} \ d\alpha \int \frac{\Vert f(x+\alpha)+f(x-\alpha)-2f(x)\Vert^{2}_{L^{4}}}{\vert \alpha \vert^{4}} \ d\alpha \right)^{1/2} \\
   & \lesssim&\Vert f \Vert_{\dot H^{2}} \Vert f_{x} \Vert_{\dot B^{1/2}_{4,2}}\Vert f \Vert_{\dot B^{3/2}_{4,2}} \\
   & \lesssim&\Vert f \Vert_{\dot H^{2}} \Vert f_{x} \Vert_{\dot H^{3/4}}\Vert f \Vert_{\dot H^{7/4}} \\
   &\lesssim&\Vert f \Vert^{2}_{\dot H^{2}} \Vert f \Vert_{\dot H^{3/2}}
   \end{eqnarray*}
   Therefore,
   \begin{equation} \label{eq:croisés2}
    \vert I_{1,3,3} \vert \lesssim \Vert f \Vert^{2}_{\dot H^{2}} \Vert f \Vert_{\dot H^{3/2}}
   \end{equation}
   
       \noindent{\bf{5.1.3.4. Estimates of $I_{1,3,4}$}} \\

  It remains to estimate $I_{1,3,4}$, to this end, we shall use the following decomposition in terms of nice controlled commutators. One may indeed rewrite $I_{1,3,4}$,  using the anti-symmetry of $\mathcal{H}$, as follows
\begin{eqnarray*}
I_{1,3,4} &=& - \int \mathcal{H}f_{xx} \ \int \frac{f_{xx}(x)}{\alpha} \  \ \int_{0}^{\infty}   e^{-\delta}  
\sin(\frac{\delta}{2}(\Delta_{\alpha}f-\bar\Delta_{\alpha}f)) \sin(\frac{\delta}{2}(\Delta_{\alpha}f+\bar\Delta_{\alpha}f))  \ d\delta \ d\alpha \ dx \\
&=& \frac{1}{2} \int f_{xx} \ \int_{0}^{\infty} \int    e^{-\delta}  \frac{1}{\alpha}
 \left[\mathcal{H} , \sin(\frac{\delta}{2}(\Delta_{\alpha}f-\bar\Delta_{\alpha}f)) \sin(\frac{\delta}{2}(\Delta_{\alpha}f+\bar\Delta_{\alpha}f)) \right]f_{xx}   \ d\delta \ d\alpha \ dx. 
     \end{eqnarray*}
 Hence, we may rewrite this term as a sum of two commutators, namely
 \begin{eqnarray*}
 I_{1,3,4} &=&  \frac{1}{2}\int\int \frac{f_{xx}(x)-f_{xx}(x-\alpha)}{\alpha} \int_{0}^{\infty}     e^{-\delta}  
 \left[\mathcal{H} , \sin(\frac{\delta}{2}(\Delta_{\alpha}f-\bar\Delta_{\alpha}f)) \sin(\frac{\delta}{2}(\Delta_{\alpha}f+\bar\Delta_{\alpha}f)) \right]f_{xx} \\
 &&\hspace{11cm} \times \ d\delta \ d\alpha \ dx \\
  &+&\frac{1}{2} \int \int \frac{f_{xx}(x-\alpha)}{\alpha} \ \int_{0}^{\infty}     e^{-\delta}  
 \left[\mathcal{H} , \sin(\frac{\delta}{2}(\Delta_{\alpha}f-\bar\Delta_{\alpha}f)) \sin(\frac{\delta}{2}(\Delta_{\alpha}f+\bar\Delta_{\alpha}f)) \right]f_{xx}   \ d\delta \ d\alpha \ dx \\
 &=&  \frac{1}{2}\int\int \frac{f_{xx}(x)-f_{xx}(x-\alpha)}{\alpha} \int_{0}^{\infty}     e^{-\delta}  
 \left[\mathcal{H} , \sin(\frac{\delta}{2}(\Delta_{\alpha}f-\bar\Delta_{\alpha}f)) \sin(\frac{\delta}{2}(\Delta_{\alpha}f+\bar\Delta_{\alpha}f)) \right]f_{xx} \\
 &&\hspace{11cm} \times \ d\delta \ d\alpha \ dx \\
  &-&\frac{1}{2} \int \int \frac{\partial_{\alpha}(f_{x}(x-\alpha)-f_x(x))}{\alpha} \ \int_{0}^{\infty}     e^{-\delta}  
 \left[\mathcal{H} , \sin(\frac{\delta}{2}(\Delta_{\alpha}f-\bar\Delta_{\alpha}f)) \sin(\frac{\delta}{2}(\Delta_{\alpha}f+\bar\Delta_{\alpha}f)) \right]f_{xx}  \\
&&\hspace{11cm} \times \ d\delta \ d\alpha \ dx.
   \end{eqnarray*}
   Finally, by integrating by parts with respect to $x$ in the first integral and with respect to $\alpha$ in the last integral, one finds,
      \begin{eqnarray} \label{I134}
I_{1,3,4} &=&- \frac{1}{2}  \int \int \frac{f_{x}(x)-f_{x}(x-\alpha)}{\alpha}  \int_{0}^{\infty}     e^{-\delta}  
 \partial_{x}\left[\mathcal{H} , \sin(\frac{\delta}{2}(\Delta_{\alpha}f-\bar\Delta_{\alpha}f)) \sin(\frac{\delta}{2}(\Delta_{\alpha}f+\bar\Delta_{\alpha}f)) \right]f_{xx} \nonumber\\
  && \hspace{11cm} \times \  d\delta \ d\alpha \ dx \nonumber  \\
  &-&\frac{1}{2} \int \int \frac{f_{x}(x-\alpha)-f_{x}(x)}{\alpha^2} \ \int_{0}^{\infty}     e^{-\delta}  
 \left[\mathcal{H} , \sin(\frac{\delta}{2}(\Delta_{\alpha}f-\bar\Delta_{\alpha}f)) \sin(\frac{\delta}{2}(\Delta_{\alpha}f+\bar\Delta_{\alpha}f)) \right]f_{xx} \nonumber \\
 && \hspace{11cm} \times \  d\delta \ d\alpha \ dx \nonumber  \\
 &+& \frac{1}{2}\int\int \frac{f_{x}(x-\alpha)-f_{x}(x)}{\alpha} \ \int_{0}^{\infty}   e^{-\delta}  
 \partial_{\alpha}\left[\mathcal{H} , \sin(\frac{\delta}{2}(\Delta_{\alpha}f-\bar\Delta_{\alpha}f)) \sin(\frac{\delta}{2}(\Delta_{\alpha}f+\bar\Delta_{\alpha}f)) \right]f_{xx} \nonumber \\
 && \hspace{11cm} \times \  d\delta \ d\alpha \ dx  \nonumber \\
 &=& I_{1,3,4,1}+I_{1,3,4,2}+I_{1,3,4,3}.
   \end{eqnarray}
   Let us estimate $I_{1,3,4,1}$. 
\begin{eqnarray*}
 \vert I_{1,3,4,1} \vert &\leq&\frac{1}{2}\int \frac{\Vert f_{x}(x)-f_{x}(x-\alpha) \Vert_{L^{2}}}{\vert \alpha \vert} \\
 &\times&\int_{0}^{\infty} e^{-\delta} 
\left\Vert  \partial_{x}\left[\mathcal{H} , \sin(\frac{\delta}{2}(\Delta_{\alpha}f-\bar\Delta_{\alpha}f)) \sin(\frac{\delta}{2}(\Delta_{\alpha}f+\bar\Delta_{\alpha}f)) \right]f_{xx} \right\Vert_{L^{\infty}} \ d \delta \ d\alpha \ dx.
\end{eqnarray*}
We then use the commutator estimate \eqref{cz} in the case $l=1$ and $k=0$. To do so, we first notice that, by using  $\vert \sin(\frac{\delta}{2}(\Delta_{\alpha}f_x + \bar \Delta_{\alpha}f_x)) \vert \leq \frac{\delta}{2} \vert \Delta_{\alpha}f_x + \bar \Delta_{\alpha}f_x \vert$ and that $\cos$ and $\sin$ functions are bounded by 1, we have for any $\delta\geq0$
\begin{eqnarray*}
 \left\Vert \partial_{x}\left(\sin(\frac{\delta}{2}(\Delta_{\alpha}f-\bar\Delta_{\alpha}f)) \sin(\frac{\delta}{2}(\Delta_{\alpha}f+\bar\Delta_{\alpha}f))\right) \right\Vert_{L^{\infty}} &\leq& \frac{\delta^2}{4} \Vert \Delta_{\alpha}f_{x}-\bar\Delta_{\alpha}f_{x} \Vert_{L^{\infty}} \Vert \Delta_{\alpha}f+\bar\Delta_{\alpha}f \Vert_{L^{\infty}} \\
 && + \ \frac{\delta}{2}\Vert \Delta_{\alpha}f_x+\bar\Delta_{\alpha}f_x \Vert_{L^{\infty}} \\
 &\leq& \frac{\delta^2}{4} \left(\Vert \Delta_{\alpha}f_x \Vert_{L^{\infty}}+\Vert \bar\Delta_{\alpha}f_x \Vert_{L^{\infty}} \right) \Vert \Delta_{\alpha}f+\bar\Delta_{\alpha}f \Vert_{L^{\infty}} \\ 
 && + \ \frac{\delta}{2}  \Vert \Delta_{\alpha}f_x+\bar\Delta_{\alpha}f_x \Vert_{L^{\infty}}. \\
\end{eqnarray*}
 Then, we find,
\begin{eqnarray*}
 \vert I_{1,3,4,1} \vert 
 &\lesssim&  \Vert f\Vert_{\dot H^{2}} \int  \frac{\Vert f_{x}(x)-f_{x}(x-\alpha) \Vert_{L^{2}}}{\vert \alpha \vert} \frac{\Vert  f_{x}(x)-f_{x}(x-\alpha)\Vert_{L^{\infty}}}{\vert \alpha \vert}   \frac{\Vert f (x-\alpha)+f(x+\alpha)-2f(x) \Vert_{L^{\infty}}}{\vert \alpha \vert }  \ d\alpha \\
 &+&\Vert f\Vert_{\dot H^{2}} \int  \frac{\Vert f_{x}(x)-f_{x}(x-\alpha) \Vert_{L^{2}}}{\vert \alpha \vert}\frac{\Vert  f_{x}(x)-f_{x}(x+\alpha)\Vert_{L^{\infty}}}{
 \vert \alpha \vert }     \frac{\Vert f (x-\alpha)+f(x+\alpha)-2f(x) \Vert_{L^{\infty}}}{\vert \alpha\vert}  \ d\alpha  \\
 &+&  \Vert f\Vert_{\dot H^{2}} \int  \frac{\Vert f_{x}(x)-f_{x}(x-\alpha) \Vert_{L^{2}}}{\vert \alpha \vert}  \frac{\Vert f_x (x-\alpha)+f_x(x+\alpha)-2f_x(x) \Vert_{L^{\infty}}}{\vert \alpha \vert }  \ d\alpha \\
  &\lesssim&  \Vert f\Vert_{\dot H^{2}}   \left(\int \frac{\Vert f_{x}(x)-f_{x}(x-\alpha) \Vert^{2}_{L^{2}}}{\vert\alpha\vert^2} \ d\alpha \right)^{1/2} \left(\int \frac{\Vert f (x-\alpha)+f(x+\alpha)-2f(x) \Vert^{4}_{L^{\infty}}}{\vert\alpha\vert^{5}}  \ d\alpha \right)^{1/4}  \\
 && \ \times \left[\left(\int\frac{\Vert  f_{x}(x)-f_{x}(x-\alpha)\Vert^{4}_{L^{\infty}}}{\vert\alpha \vert^3} \ d\alpha \right)^{1/4}+ \left(\int \frac{\Vert  f_{x}(x)-f_{x}(x+\alpha)  \Vert^{4}_{L^{\infty}}}{\vert\alpha\vert^3} \ d\alpha \right)^{1/4} \right]  \\
 && + \ \Vert f\Vert_{\dot H^{2}}   \left(\int \frac{\Vert f_{x}(x)-f_{x}(x-\alpha) \Vert^{2}_{L^{2}}}{\vert\alpha\vert^2} \ d\alpha \right)^{1/2} \left(\int \frac{\Vert f (x-\alpha)+f(x+\alpha)-2f(x) \Vert^{2}_{L^{\infty}}}{\vert\alpha\vert^{2}}  \ d\alpha \right)^{1/2}  \\
 &\lesssim&  \Vert f\Vert_{\dot H^{2}}
  \Vert f \Vert_{\dot H^{3/2}} \Vert f \Vert_{\dot B^{1}_{\infty,4}} \Vert f_{x} \Vert_{\dot B^{1/2}_{\infty,4}} + \Vert f\Vert_{\dot H^{2}}
  \Vert f \Vert_{\dot H^{3/2}} \Vert f \Vert_{\dot B^{1}_{\infty,2}} \\
 &\lesssim&  \Vert f\Vert^{2}_{\dot H^{2}} \Vert f \Vert^{2}_{\dot H^{3/2}},
 \end{eqnarray*}
 where we used the following embeddings $\dot H^{3/2} \hookrightarrow \dot B^{1}_{\infty,4}$, $\dot H^{2} \hookrightarrow \dot B^{3/2}_{\infty,4}$,  and $\dot H^{3/2} \hookrightarrow \dot B^{1}_{\infty,2}$. \\
 
Then, we estimate $I_{1,3,4,2}$, we see that this term  may be written as
\begin{eqnarray} \label{I1342}
I_{1,3,4,2}&=& \frac{1}{2}\int \mathcal{H}f_{xx} \frac{f_{x}(x-\alpha)-f_{x}(x)}{\alpha^2}  \int_{0}^{\infty} \int    e^{-\delta}  \sin(\frac{\delta}{2}(\Delta_{\alpha}f-\bar\Delta_{\alpha}f)) \sin(\frac{\delta}{2}(\Delta_{\alpha}f+\bar\Delta_{\alpha}f)) \  d\delta  \ d\alpha \ dx \nonumber \\
&+&\frac{1}{2} \int f_{xx} \frac{\mathcal{H}f_{x}(x-\alpha)-\mathcal{H}f_{x}(x)}{\alpha^2}  \int_{0}^{\infty} \int    e^{-\delta}  \sin(\frac{\delta}{2}(\Delta_{\alpha}f-\bar\Delta_{\alpha}f)) \sin(\frac{\delta}{2}(\Delta_{\alpha}f+\bar\Delta_{\alpha}f))  \  d\delta \ d\alpha \ dx \nonumber\\
&=& I_{1,3,4,2,1}+ I_{1,3,4,2,2}.
\end{eqnarray}
To estimate $I_{1,3,4,2,1}$, we write
\begin{eqnarray*}
\vert I_{1,3,4,2,1} \vert &\leq&\frac{1}{2} \Vert f \Vert_{\dot H^{2}} \int \frac{\Vert f_{x}(x-\alpha)-f_{x}(x) \Vert_{L^{2}}}{\alpha^2} \ \int_{0}^{\infty}     e^{-\delta}  \frac{\Vert \delta_{\alpha}f+\bar\delta_{\alpha}f \Vert_{L^{\infty}}}{\vert \alpha \vert}  \  \ d\delta \ d\alpha   \\
& \lesssim& \Vert f \Vert_{\dot H^{2}} \left(\int \frac{\Vert f_{x}(x-\alpha)-f_{x}(x) \Vert^{2}_{L^{2}}}{\alpha^2} \ d\alpha \right)^{1/2}  \left( \int \frac{\Vert \delta_{\alpha}f+\bar\delta_{\alpha}f \Vert^{2}_{L^{\infty}}}{\alpha^4} \ d\alpha \right)^{1/2} \\
& \lesssim& \Vert f \Vert_{\dot H^{2}} \Vert f \Vert_{\dot H^{3/2}}  \Vert f \Vert_{\dot B^{3/2}_{\infty,2}} \\
& \lesssim& \Vert f \Vert^{2}_{\dot H^{2}} \Vert f \Vert_{\dot H^{3/2}}  
\end{eqnarray*}
Since $\mathcal{H}$ is continuous on $L^{2}$, one gets the same control for $I_{1,3,4,2,2}$ so that one finally obtains
\begin{equation*}
\vert I_{1,3,4,2} \vert \lesssim \Vert f \Vert^{2}_{\dot H^{2}} \Vert f \Vert_{\dot H^{3/2}}
\end{equation*}
It therefore remains to estimate $I_{1,3,4,3}$. To do so, we open the commutator and use the anti-symmetry of $\mathcal{H}$. Then, we explicitely compute the derivative in $\alpha$. Hence, we obtain the following decomposition
\begin{eqnarray} \label{I1343}
 I_{1,3,4,3}&=&-\frac{1}{4}\int f_{xx} \frac{\mathcal{H}f_{x}(x-\alpha)-\mathcal{H}f_{x}(x)}{\alpha} \ \int_{0}^{\infty} \int  \delta  e^{-\delta}  \partial_{\alpha}(\Delta_{\alpha}f-\bar\Delta_{\alpha}f) \ \cos(\frac{\delta}{2}(\Delta_{\alpha}f-\bar\Delta_{\alpha}f)) \nonumber \\
 && \times \sin(\frac{\delta}{2}(\Delta_{\alpha}f+\bar\Delta_{\alpha}f))   \ d\delta \ d\alpha \ dx \nonumber \\
 &-&\frac{1}{4}\int \mathcal{H}f_{xx} \frac{f_{x}(x-\alpha)-f_{x}(x)}{\alpha} \ \int_{0}^{\infty} \int  \delta  e^{-\delta}  \partial_{\alpha}(\Delta_{\alpha}f-\bar\Delta_{\alpha}f) \ \cos(\frac{\delta}{2}(\Delta_{\alpha}f-\bar\Delta_{\alpha}f)) \\
 && \times \sin(\frac{\delta}{2}(\Delta_{\alpha}f+\bar\Delta_{\alpha}f)    \ d\delta \ d\alpha \ dx \nonumber \\
 &-&\frac{1}{4} \int f_{xx} \int \frac{\mathcal{H}f_{x}(x-\alpha)-\mathcal{H}f_{x}(x)}{\alpha} \ \int_{0}^{\infty}   \delta  e^{-\delta}  \partial_{\alpha}(\Delta_{\alpha}f+\bar\Delta_{\alpha}f) \ \cos(\frac{\delta}{2}(\Delta_{\alpha}f+\bar\Delta_{\alpha}f))  \nonumber\\
 && \times \sin(\frac{\delta}{2}(\Delta_{\alpha}f-\bar\Delta_{\alpha}f)  \  \ d\delta \ d\alpha \ dx \nonumber\\
 &-&\frac{1}{4} \int \mathcal{H}f_{xx} \int \frac{f_{x}(x-\alpha)-f_{x}(x)}{\alpha} \ \int_{0}^{\infty}   \delta  e^{-\delta}  \partial_{\alpha}(\Delta_{\alpha}f+\bar\Delta_{\alpha}f) \ \cos(\frac{\delta}{2}(\Delta_{\alpha}f+\bar\Delta_{\alpha}f)) \nonumber\\
 && \times \sin(\frac{\delta}{2}(\Delta_{\alpha}f-\bar\Delta_{\alpha}f)  \  \ d\delta \ d\alpha \ dx \nonumber \\
 &=& \sum_{i=1}^{4}  I_{1,3,4,3,i}.
\end{eqnarray}
Since we shall do $L^{p}$ estimates, $p\in(1,\infty)$ on the terms involving $\mathcal{H}$ one observes that these terms have the same regularity as the terms $I_{1,3,j}$ for $j=2,3$ (see \eqref{I132}). It is therefore not difficult to see that,
$$
\left \vert I_{1,3,4,3,1}+I_{1,3,4,3,2} \right\vert \lesssim  \Vert f \Vert^{2}_{\dot H^{2}} \Vert f \Vert^{2}_{\dot H^{3/2}},
$$
and,
$$
\left \vert I_{1,3,4,3,3}+I_{1,3,4,3,4} \right\vert \lesssim \Vert f \Vert^{2}_{\dot H^{2}} \Vert f \Vert_{\dot H^{3/2}}. \\
$$
Therefore, we have obtained
\begin{equation} \label{i13}
\vert I_{1,3} \vert \lesssim \Vert f \Vert^{2}_{\dot H^{2}} \left(\Vert f \Vert^{2}_{\dot H^{3/2}} +  \Vert f \Vert_{\dot H^{3/2}} \right).
\end{equation}

   \noindent {\bf{5.2. Estimate of $I_2$}} \\
  
  To estimate
  \begin{eqnarray*}
 I_{2}&=&  -\int  \mathcal{H}f_{xx} \ \int(\partial_{x} \Delta_{\alpha} f)^{2}  \  \ \int_{0}^{\infty} \delta  e^{-\delta}  \sin(\delta\Delta_{\alpha}f(x) ) \ d\delta \ d\alpha \ dx,   
  \end{eqnarray*}
  it suffices to write that,
\begin{equation*}
  \vert I_{2} \vert \leq \Gamma(2)\Vert f\Vert_{\dot H^{2}} \left(\int \frac{\Vert f_{x}(x)-f_{x}(x-\alpha) \Vert^{2}_{L^{2}}}{\alpha^{2}} \ d\alpha\right)^{1/2} \left(\int \frac{\Vert f_{x}(x)-f_{x}(x-\alpha) \Vert^{2}_{L^{\infty}}}{\alpha^{2}} \ d\alpha\right)^{1/2}.
\end{equation*}
Therefore, we obtain
  \begin{eqnarray} \label{est2}
  \vert I_{2} \vert & \lesssim& \Vert f \Vert_{\dot H^{2}} \Vert f_{x}\Vert_{\dot B^{1/2}_{\infty,2}} \Vert f \Vert_{\dot H^{3/2}} \nonumber \\
  & \lesssim& \Vert f \Vert^{2}_{\dot H^{2}} \Vert f \Vert_{\dot H^{3/2}}
  \end{eqnarray}

\noindent Finally, combining all the estimates, we get that
  \begin{equation}
 \frac{1}{2} \partial_{t} \Vert f \Vert^{2}_{\dot H^{3/2}} + \frac{\pi}{1+K^2} \Vert f \Vert^{2}_{\dot H^{2}} \lesssim  \Vert f \Vert^{2}_{\dot H^{2}} \left(\Vert f \Vert^{2}_{\dot H^{3/2}} + \Vert f \Vert_{ \dot H^{3/2}}\right).
  \end{equation}
  And then integrating in time $s \in [0,T]$ one gets the desired energy inequality of Lemma \ref{h3/2}.
  Therefore, if $\Vert f_{0} \Vert_{\dot H^{3/2}}$ is smaller than some $C(K)$ that depends only on $K$, then the solution is in $L^{\infty}([0,T],L^{2}) \cap L^{2}([0,T], \dot H^{2})$.  This concludes the $\dot H^{3/2}$-estimates and therefore Lemma \ref{h3/2} is proved.
   \qed
  
  \section{A priori estimates in $\dot H^{5/2}$ estimates} \label{h52}
  
  In this section we shall prove the following lemma.
\begin{lemma} \label{ineg2}
Let $T>0$ and $f_{0} \in \dot H^{5/2} \cap  \dot H^{3/2}$ so that $\Vert f_{0} \Vert_{\dot H^{3/2}} < C(\Vert \partial_{x}f_{0} \Vert_{L^{\infty}})$, then we have 
\begin{equation*} \label{ineg}
 \Vert f \Vert^{2}_{\dot H^{5/2}}(T) + \frac{\pi}{1+K^2} \int_{0}^T \Vert f \Vert^{2}_{\dot H^{3}} \ ds \lesssim \Vert f_{0} \Vert_{\dot H^{5/2}} +  \left(\Vert f \Vert_{L^{\infty}([0,T],\dot H^{3/2})} + \Vert f \Vert^{2}_{L^{\infty}([0,T],\dot H^{3/2})}  \right) \int_{0}^T \Vert f \Vert^{2}_{\dot H^{3}} \ ds
\end{equation*}
where $K=\displaystyle\sup_{(x,t)\in \mathbb R \times  [0,T]}\vert f_{x}(x,t)\vert$.
\end{lemma}
\noindent {\bf{Proof of Lemma \ref{ineg2}}} We now do {\it{a priori}} estimate in $\dot H^{5/2}$. Using the anti-symmetry property of the Hilbert transform, one finds

\begin{eqnarray*}
 \frac{1}{2}\partial_{t} \Vert f \Vert^{2}_{\dot H^{5/2}}  
 &=&   \int \mathcal{H}\partial^{3}_{x} f 
\ \int \partial^{3}_{x} \Delta_{\alpha}f \int_{0}^{\infty} e^{-\delta}  \ \cos(\delta \Delta_{\alpha} f)  \ dx \ d\alpha \ d\delta \\
  &-&3\int \mathcal{H}\partial^{3}_{x} f 
\ \int \partial_{x} \Delta_{\alpha}f \ \partial_{xx} \Delta_{\alpha}f \int_{0}^{\infty} \delta e^{-\delta}  \ \sin(\delta \Delta_{\alpha} f)  \ dx \ d\alpha \ d\delta \\
&-& \int \mathcal{H}\partial^{3}_{x} f 
\ \int (\partial_{x} \Delta_{\alpha}f)^3 \int_{0}^{\infty} \delta^{2} e^{-\delta}  \ \cos(\delta \Delta_{\alpha} f)  \ dx \ d\alpha \ d\delta \\
&=& T_1 + T_2 + T_3.
  \end{eqnarray*}
  Analogously to the $\dot H^{3/2}$ {\it{a priori}} estimate, we decompose the first term as follows. By noticing that $I_1$ is analogous to $T_1$, we find that
  \begin{eqnarray*}
T_1&=& -\int \mathcal{H}\partial^{3}_{x}f \ \int (\partial^{3}_{x}  \Delta_{\alpha} f -\partial^{3}_{x} \bar \Delta_{\alpha} f)  \  \ \int_{0}^{\infty}   e^{-\delta} \cos(\frac{\delta}{2}(\Delta_{\alpha}f-\bar\Delta_{\alpha}f)) \sin^{2}(\frac{\delta}{4}(\Delta_{\alpha}f+\bar\Delta_{\alpha}f)) \ d\delta \ d\alpha \ dx \\
&+&\frac{1}{2} \int \mathcal{H}\partial^{3}_{x}f \ \int (\partial^{3}_{x}  \Delta_{\alpha} f -\partial^{3}_{x} \bar \Delta_{\alpha} f)  \  \ \int_{0}^{\infty}   e^{-\delta} \cos(\frac{\delta}{2}(\Delta_{\alpha}f-\bar\Delta_{\alpha}f))\ d\delta \ d\alpha \ dx \\
&+& \int \mathcal{H}\partial^{3}_{x}f \ \int \frac{\partial^{3}_{x}f(x)-\partial^{3}_{x}f(x-\alpha)}{\alpha} \  \ \int_{0}^{\infty}   e^{-\delta}  
\sin(\frac{\delta}{2}(\Delta_{\alpha}f-\bar\Delta_{\alpha}f)) \sin(\frac{\delta}{2}(\Delta_{\alpha}f+\bar\Delta_{\alpha}f))  \ d\delta \ d\alpha \ dx \\
&=& T_{1,1}+T_{1,2}+T_{1,3}.
\end{eqnarray*}
Note that $T_{1,1}$, $T_{1,2}$, and $T_{1,3}$ are respectively analogous to $I_{1,1}$, $I_{1,2}$ and $I_{1,3}$. Using \eqref{I13}, we immedialy infer that
  \begin{eqnarray*}
T_1&=& -\int \mathcal{H}\partial^{3}_{x}f \ \int (\partial^{3}_{x}  \Delta_{\alpha} f -\partial^{3}_{x} \bar \Delta_{\alpha} f)  \  \ \int_{0}^{\infty}   e^{-\delta} \cos(\frac{\delta}{2}(\Delta_{\alpha}f-\bar\Delta_{\alpha}f)) \sin^{2}(\frac{\delta}{4}(\Delta_{\alpha}f+\bar\Delta_{\alpha}f)) \ d\delta \ d\alpha \ dx \\
&+&\frac{1}{2} \int \mathcal{H}\partial^{3}_{x}f \ \int (\partial^{3}_{x}  \Delta_{\alpha} f -\partial^{3}_{x} \bar \Delta_{\alpha} f)  \  \ \int_{0}^{\infty}   e^{-\delta} \cos(\frac{\delta}{2}(\Delta_{\alpha}f-\bar\Delta_{\alpha}f))\ d\delta \ d\alpha \ dx \\
&+& \int \mathcal{H}\partial^{3}_{x}f \ \int \frac{f_{xx}(x)-f_{xx}(x-\alpha)}{\alpha^2} \  \ \int_{0}^{\infty}   e^{-\delta}  
\sin(\frac{\delta}{2}(\Delta_{\alpha}f-\bar\Delta_{\alpha}f)) \sin(\frac{\delta}{2}(\Delta_{\alpha}f+\bar\Delta_{\alpha}f))  \ d\delta \ d\alpha \ dx \\
&-& \frac{1}{2} \int \mathcal{H}\partial^{3}_{x}f \ \int \frac{f_{xx}(x)-f_{xx}(x-\alpha)}{\alpha} \  \ \int_{0}^{\infty} \delta  e^{-\delta}  
\partial_{\alpha} D \ \cos(\frac{\delta}{2}(\Delta_{\alpha}f-\bar\Delta_{\alpha}f)) \sin(\frac{\delta}{2}(\Delta_{\alpha}f+\bar\Delta_{\alpha}f))  \\
&& \hspace{13cm} \times \ d\delta \ d\alpha \ dx \\
&-& \frac{1}{2} \int \mathcal{H}\partial^{3}_{x}f \ \int \frac{f_{xx}(x)-f_{xx}(x-\alpha)}{\alpha} \  \ \int_{0}^{\infty} \delta  e^{-\delta}  
 \ \partial_{\alpha} S \ \sin(\frac{\delta}{2}(\Delta_{\alpha}f-\bar\Delta_{\alpha}f))  \ \cos(\frac{\delta}{2}(\Delta_{\alpha}f+\bar\Delta_{\alpha}f))  \\
 &&  \hspace{13cm} \times d\delta \ d\alpha \ dx \\
 &-& \int \mathcal{H}\partial^{3}_{x}f \ \int \partial^{3}_{x}f(x) \  \ \int_{0}^{\infty}   e^{-\delta}  
\sin(\frac{\delta}{2}(\Delta_{\alpha}f-\bar\Delta_{\alpha}f)) \sin(\frac{\delta}{2}(\Delta_{\alpha}f+\bar\Delta_{\alpha}f)) \frac{1}{\alpha} \ d\delta \ d\alpha \ dx \\
&=& \sum_{i=1}^{6} T_{1,j}
\end{eqnarray*}

We then estimate the $T_{1,j}$, $j=1,...,6$. \\

\noindent {\bf{6.1. Estimate  of  $T_{1,1}$}} \\

To control $T_{1,1}$, we use the continuity of the Hilbert transform on $L^{2}$ along with the embedding $\dot H^{3/2} \hookrightarrow \dot B^{1}_{\infty,2}$, then one gets,
\begin{eqnarray*}
T_{1,1}&=&-\int \mathcal{H}\partial^{3}_{x}f \ \int (\partial^{3}_{x}  \Delta_{\alpha} f -\partial^{3}_{x} \bar \Delta_{\alpha} f)  \  \ \int_{0}^{\infty}  e^{-\delta} \cos(\frac{\delta}{2}(\Delta_{\alpha}f-\bar\Delta_{\alpha}f)) \sin^{2}(\frac{\delta}{4}(\Delta_{\alpha}f+\bar\Delta_{\alpha}f)) \ d\delta \ d\alpha \ dx \\
&\leq& \frac{\Gamma(3)}{4}  \Vert f \Vert^{2}_{\dot H^{3}} \int \frac{\Vert \delta_{\alpha}f+\bar\delta_{\alpha}f \Vert^{2}_{L^{\infty}}}{\vert \alpha \vert^{3}} \ d\alpha \\
&\lesssim& \Vert f \Vert^{2}_{\dot H^{3}} \Vert f \Vert^{2}_{\dot B^{1}_{\infty,2}} \\
&\lesssim&  \Vert f \Vert^{2}_{\dot H^{3}} \Vert f \Vert^{2}_{\dot H^{3/2}}
 \end{eqnarray*}

\noindent {\bf{6.2. Estimate  of  $T_{1,2}$}}  \\

 We first rewrite $T_{1,2}$ as follows, by integrating by parts, one finds
 \begin{eqnarray*}
T_{1,2}&=&\frac{1}{2}\int \mathcal{H}\partial^{3}_{x}f \ \int (  \partial^{3}_{x}\Delta_{\alpha} f - \partial^{3}_{x} \bar \Delta_{\alpha} f)  \  \ \int_{0}^{\infty}   e^{-\delta} \cos(\frac{\delta}{2}(\Delta_{\alpha}f-\bar\Delta_{\alpha}f))\ d\delta \ d\alpha \ dx \\
&=&\frac{1}{2} \int \mathcal{H}\partial^{3}_{x}f \ \int \frac{1}{\alpha} \partial_{\alpha}\{\delta_{\alpha} f_{xx} + \bar \delta_{\alpha} f_{xx} \}  \  \ \int_{0}^{\infty}   e^{-\delta} \cos(\frac{\delta}{2}(\Delta_{\alpha}f-\bar\Delta_{\alpha}f))\ d\delta \ d\alpha \ dx \\
&=& \frac{1}{2} \int \mathcal{H}\partial^{3}_{x}f  \int \frac{f_{xx}(x-\alpha)+f_{xx}(x+\alpha)-2f_{xx}(x)}{\alpha^2}    \int_{0}^{\infty}   e^{-\delta} \cos(\frac{\delta}{2} (\bar\Delta_{\alpha}f-\Delta_{\alpha}f))\ d\delta \ d\alpha \ dx \\
&+& \frac{1}{4}\int \mathcal{H}\partial^{3}_{x}f  \int \frac{f_{xx}(x-\alpha)+f_{xx}(x+\alpha)-2f_{xx}(x)}{\alpha}     \int_{0}^{\infty} \delta  e^{-\delta} \\
&&\times \ \partial_{\alpha}\{ \Delta_{\alpha}f-\bar\Delta_{\alpha}f\} \sin(\frac{\delta}{2} (\Delta_{\alpha}f-\bar\Delta_{\alpha}f))\ d\delta \ d\alpha \ dx. \\
\end{eqnarray*}

\noindent Then, we obtain that,
\begin{eqnarray} \label{dissip}
T_{1,2}&=&\nonumber - \int \mathcal{H}\partial^{3}_{x}f \ \int \frac{f_{xx}(x-\alpha)+f_{xx}(x+\alpha)-2f_{xx}(x)}{\alpha^2}    \  \ \int_{0}^{\infty}   e^{-\delta} \sin^{2}(\frac{\delta}{4} (\bar\Delta_{\alpha}f-\Delta_{\alpha}f))\ d\delta \ d\alpha \ dx \\
&+& \frac{1}{4}\nonumber \int \mathcal{H}\partial^{3}_{x}f \ \int \frac{f_{xx}(x-\alpha)+f_{xx}(x+\alpha)-2f_{xx}(x)}{\alpha}    \  \ \int_{0}^{\infty} \delta  e^{-\delta} \frac{f_{x}(x-\alpha)+f_{x}(x+\alpha)-2f_{x}(x)}{\alpha} \\
&\times&\nonumber \sin(\frac{\delta}{2} (\Delta_{\alpha}f-\bar\Delta_{\alpha}f))\ d\delta \ d\alpha \ dx \\
&-&\frac{1}{4} \nonumber \int \mathcal{H}\partial^{3}_{x}f \ \int \frac{f_{xx}(x-\alpha)+f_{xx}(x+\alpha)-2f_{xx}(x)}{\alpha^3}    \  \ \int_{0}^{\infty} \delta  e^{-\delta} \int_{0}^{\alpha} f_{x}(x-s)+f_{x}(x+s)-2f_{x}(x) \ ds \\
&\times&\nonumber \sin(\frac{\delta}{2} (\Delta_{\alpha}f-\bar\Delta_{\alpha}f))\ d\delta \ d\alpha \ dx \\
&-&  \pi\int \mathcal{H}\partial^{3}_{x}f \Lambda f_{xx} \ dx \\
&=&\nonumber T_{1,2,1}+T_{1,2,2}+T_{1,2,3}+T_{1,2,4}
 \end{eqnarray} 
 In order to estimate $T_{1,2}$ we shall estimate the $T_{1,2,j}$, for $j=1,2,3$. Note that $T_{1,2,4}$ is a dissipative term. We start with  $T_{1,2,1}$, we observe that this term is analogous to $I_{1,2,1}$ (see \eqref{10sip}). More precisely, we write

\begin{eqnarray*}
T_{1,2,1}&=&- \int \mathcal{H}\partial^{3}_{x}f \ 
\int \frac{f_{xx}(x-\alpha)+f_{xx}(x+\alpha)-2f_{xx}(x)}{\alpha^2}    \  \ \int_{0}^{\infty}   e^{-\delta} 
\sin^{2}(\frac{\delta}{4} (\bar\Delta_{\alpha}f-\Delta_{\alpha}f))\ d\delta \ d\alpha \ dx  \\
&=&- \int \mathcal{H}\partial^{3}_{x}f \ \int \frac{f_{xx}(x-\alpha)+f_{xx}(x+\alpha)-2f_{xx}(x)}{\alpha^2}    \  \\
&\times&  \int_{0}^{\infty}   e^{-\delta} \left(\sin^{2}(\frac{\delta}{4} (\bar\Delta_{\alpha}f-\Delta_{\alpha}f))-\sin^{2}(\frac{\delta}{2} f_{x}(x))\right)\ d\delta \ d\alpha \ dx \\
&+&  \int \mathcal{H}\partial^{3}_{x}f \ \int \frac{f_{xx}(x-\alpha)+f_{xx}(x+\alpha)-2f_{xx}(x)}{\alpha^2}  \int_{0}^{\infty}   e^{-\delta}  \sin^{2}(\frac{\delta}{2} f_{x}(x))\ d\delta \ d\alpha \ dx \\
&=&T_{1,2,1,1}+T_{1,2,1,2}.
 \end{eqnarray*} 
In order to estimate $T_{1,2,1,1}$ we use the estimate $\eqref{inegtrigo}$ together with the identity \eqref{formuleD}. So that by using  H\"older inequality  one finds that
  \begin{eqnarray*}
 \vert T_{1,2,1,1}\vert &\leq&  \frac{\Gamma(2)}{4} \Vert f \Vert_{\dot H^{3}} \int \frac{\Vert f_{xx}(x-\alpha)+f_{xx}(x+\alpha)-2f_{xx}(x)\Vert_{L^{\infty}}}{\vert \alpha \vert^{3}}  \int_{0}^{\alpha}  {\Vert f_{x}(x-s)-f_{x}(x) \Vert_{L^{2}}} \ ds \ d\alpha.  \\
 \end{eqnarray*}
 Then, by using H\"older inequality in the $s$-integral (where $r$ and $\bar r$ are conjugate) and for some real number $q \in [0,2]$ that will be chosen later, we have
 \begin{eqnarray*}
  \vert T_{1,2,1,1}\vert  &\lesssim&   \Vert f \Vert_{\dot H^{3}} \int \frac{\Vert f_{xx}(x-\alpha)-f_{xx}(x) \Vert_{L^{\infty}}}{\vert \alpha \vert^{3}} \vert \alpha \vert^{q+\frac{1}{\bar r}}   \left(\int  \frac{\Vert f_{x}(x-s)-f_{x}(x) \Vert^{r} _{L^{2}}}{\vert s \vert^{qr}} \ ds\right)^{1/r}  \ d\alpha  \\
 &\lesssim&   \Vert f \Vert_{\dot H^{3}} \Vert f_{xx} \Vert_{\dot B^{2-q-\frac{1}{\bar r}}_{\infty,1}} \Vert f_{x} \Vert_{\dot B^{q-\frac{1}{r}}_{2,r}}
  \end{eqnarray*}
 By choosing $q={9/8},$  $r=\bar r=2,$ we find
  \begin{eqnarray*}
  \vert T_{1,2,1,1} \vert &\leq&   \Vert f \Vert_{\dot H^{3}} \Vert f_{xx} \Vert_{\dot B^{3/8}_{\infty,1}} \Vert f_{x} \Vert_{\dot B^{5/8}_{2,2}}
 \end{eqnarray*}
  Then, by  interpolating in the space  $\dot B^{19/8}_{\infty,1}$, we find,
 \begin{eqnarray*}
 \vert T_{1,2,1,1} \vert &\leq&  \Vert f \Vert_{\dot H^{3}} \Vert f\Vert_{\dot B^{13/8}_{2,2}}  \Vert f \Vert^{1/12}_{\dot B^{1}_{\infty,\infty}}  \Vert f \Vert^{11/12}_{\dot B^{5/2}_{\infty,\infty}}  \\
 \end{eqnarray*}
 Since,
 $$
  \Vert f\Vert_{\dot B^{13/8}_{2,2}} \leq \Vert f \Vert^{1/12}_{\dot H^{3}} \Vert f \Vert^{11/12}_{\dot H^{3/2}},
  $$
then, by using classical Besov embeddings, one finally gets
$$
\vert T_{1,2,1,1} \vert \leq  \Vert f \Vert^{2}_{\dot H^{3}}  \Vert f \Vert_{\dot H^{3/2}} 
$$
As for the term $T_{1,2,1,2}$, using the same idea as the estimate of the term $I_{1,2,1,2}$ (see \eqref{I1212}). We have that,
\begin{eqnarray} \label{+}
 T_{1,2,1,2} &=&\nonumber -  \int \mathcal{H}\partial^{3}_{x}f \ \int \frac{f_{xx}(x-\alpha)+f_{xx}(x+\alpha)-2f_{xx}(x)}{\alpha^2}    \  \ \int_{0}^{\infty}   e^{-\delta} \ \sin^{2}( \frac{\delta}{2} f_{x}(x)) \ d\delta \ d\alpha \ dx   \\
&=& 2\pi   \int (\mathcal{H}\partial^{3}_{x}f)^2     \ \int_{0}^{\infty}   e^{-\delta} \ \sin^{2}( \frac{\delta}{2} f_{x}(x)) \ d\delta \  \ dx   \leq \pi \frac{K^2}{1+ K^2} \Vert f \Vert^{2}_{\dot H^{3}},
\end{eqnarray}
where $K=\displaystyle\sup_{(x,t)\in \mathbb R \times  [0,T]}\vert f_{x}(x,t)\vert$. Therefore, we find that for some fixed constant $C>0$,
\begin{equation}
 T_{1,2,1,1}+T_{1,2,1,2}\leq C \Vert f \Vert^{2}_{\dot H^{3}}  \Vert f \Vert_{\dot H^{3/2}}  + \pi \frac{ K^2}{1+  K^2} \Vert f \Vert^{2}_{\dot H^{3}}.
\end{equation}

To estimate  $T_{1,2,2}$, it suffices to write that

\begin{eqnarray*}
T_{1,2,2} &=&-\int \mathcal{H}\partial^{3}_{x}f \ \int \frac{f_{xx}(x-\alpha)+f_{xx}(x+\alpha)-2f_{xx}(x)}{\alpha}    \  \ \int_{0}^{\infty} \delta  e^{-\delta} \frac{f_{x}(x-\alpha)+f_{x}(x+\alpha)-2f_{x}(x)}{\alpha} \\
&\times& \sin(\frac{\delta}{2} (\Delta_{\alpha}f-\bar\Delta_{\alpha}f))\ d\delta \ d\alpha \ dx \\
&\leq& \Vert f \Vert^{2}_{\dot H^{3}} 
\left(\int \frac{\Vert f_{xx}(x-\alpha)+f_{xx}(x+\alpha)-2f_{xx}(x)\Vert^{2}_{L^{\infty}}}{\alpha^2} \frac{ \Vert f_{x}(x-\alpha)+f_{x}(x+\alpha)-2f_{x}(x)\Vert_{L^{2}}}{\alpha^2} \ d\alpha \right)^{1/2} \\
&\leq& \Vert f \Vert_{\dot H^{3}} \Vert f \Vert_{\dot B^{5/2}_{\infty,2}}  \Vert f_{x} \Vert_{\dot B^{1/2}_{2,2}}\\
&\leq& \Vert f \Vert^{2}_{\dot H^{3}} \Vert f \Vert_{\dot H^{3/2}}  \\
\end{eqnarray*}

\noindent We now estimate $T_{1,2,3}$.  Let $q \in [0,2]$ and let $r, \bar r$ so that $1/r+1/\bar r=1$, we write

\begin{eqnarray*}
T_{1,2,3}&=&-\frac{1}{4}\int \mathcal{H}\partial^{3}_{x}f \ \int \frac{f_{xx}(x-\alpha)+f_{xx}(x+\alpha)-2f_{xx}(x)}{\alpha^3}    \  \ \int_{0}^{\infty} \delta  e^{-\delta}  \\
&& \ \times \int_{0}^{\alpha} f_{x}(x-s)+f_{x}(x+s)-2f_{x}(x)  \ ds\sin(\frac{\delta}{2} (\Delta_{\alpha}f-\bar\Delta_{\alpha}f))\ d\delta \ d\alpha \ dx \\
&\leq&\frac{1}{4}\Vert f \Vert_{\dot H^{3}} \int_{0}^{\infty} \delta e^{-\delta} \frac{\Vert f_{xx}(x-\alpha)+f_{xx}(x+\alpha)-2f_{x}(xx) \Vert_{L^{\infty}}}{\vert \alpha \vert^{3}} \\
 && \ \times \int_{0}^{\alpha}  {\Vert f_{x}(x-s)-f_{x}(x) \Vert_{L^{2}}}+{\Vert f_{x}(x+s)-f_{x}(x) \Vert_{L^{2}}} \ ds  \ d\delta \ d\alpha \\
 &\lesssim& \Vert f \Vert_{\dot H^{3}} \left( \int \frac{\Vert f_{xx}(x-\alpha)+ f_{xx}(x+\alpha)-2f_{xx}(x) \Vert_{L^{\infty}}}{\vert \alpha \vert^{3-q-\frac{1}{\bar r}}}  \ d\alpha \right)  \\
  && \ \times \left[\left(\int  \frac{\Vert f_{x}(x-s)-f_{x}(x) \Vert^{r} _{L^{2}}}{\vert s \vert^{qr}} \ ds\right)^{1/r} + \left(\int  \frac{\Vert f_{x}(x+s)-f_{x}(x) \Vert^{r} _{L^{2}}}{\vert s \vert^{qr}} \ ds\right)^{1/r}\right]   \\
   &\lesssim&  \Vert f \Vert_{H^{3}} \Vert f_{xx} \Vert_{\dot B^{2-q-\frac{1}{\bar r}}_{\infty,1}} \Vert f_{x} \Vert_{\dot B^{q-\frac{1}{r}}_{2,r}}
\end{eqnarray*}
Then, by choosing  $q={2},$ $r=\bar r=2$, we get
\begin{eqnarray*}
\vert T_{1,2,3} \vert&\lesssim&\Vert f \Vert_{H^{3}} \Vert f \Vert_{\dot B^{3/2}_{\infty,1}} \Vert f \Vert_{\dot B^{5/2}_{2,2}} 
\end{eqnarray*}

Then, using the following interpolations $\dot H^{5/2}=\left[\dot H^{3/2}, \dot H^{3}\right]_{\frac{1}{3},\frac{2}{3}}$ and  $\dot B^{3/2}_{\infty,1}=\left[\dot B^{1}_{\infty,\infty},, \dot B^{5/2}_{\infty,\infty}\right]_{\frac{2}{3},\frac{1}{3}}$, we find
\begin{eqnarray*}
\vert T_{1,2,3} \vert & \lesssim&   \Vert f \Vert_{H^{3}} \Vert f \Vert^{2/3}_{\dot H^{3}}\Vert f \Vert^{1/3}_{\dot H^{3/2}}  \Vert f \Vert^{1/3}_{\dot B^{5/2}_{\infty,\infty}} \Vert f \Vert^{2/3}_{\dot B^{1}_{\infty,\infty}} \\
& \lesssim&  \Vert f \Vert^{2}_{H^{3}} \Vert f \Vert_{\dot H^{3/2}} \\
\end{eqnarray*}

 \noindent {\bf{6.3. Estimate  of  $T_{1,3}$}}  \\
 
 To estimate $T_{1,3}$, it suffices to write that
 \begin{eqnarray*}
\vert T_{1,3} \vert &\leq&  \Vert f \Vert_{\dot H^{3}}  \  \ \int_{0}^{\infty} \delta  e^{-\delta} \frac{\Vert f_{xx}(x)- f_{xx}(x-\alpha)\Vert_{L^{\infty}}\Vert f(x-\alpha)+f(x+\alpha)-2f(x)\Vert_{L^{2}}}{\vert \alpha \vert^{3}} \ d\delta \ d\alpha  \\
&\leq&  \Gamma(2) \Vert f \Vert_{\dot H^{3}} \left(\int \frac{\Vert  f_{xx}(x)- f_{xx}(x-\alpha) \Vert^{2}_{L^{\infty}}}{\vert \alpha \vert^2} \ d\alpha \int \frac{\Vert f(x-\alpha)+f(x+\alpha)-2f(x) \Vert^{2}_{L^{2}}}{\vert \alpha \vert^4} \ d\alpha \right)^{1/2} \\
&\leq&  \Vert f \Vert_{\dot H^{3}} \Vert f_{xx} \Vert_{\dot B^{1/2}_{\infty,2}} \Vert f \Vert_{\dot B^{3/2}_{2,2}} \\
&\leq& \Vert f \Vert^{2}_{\dot H^{3}} \Vert f \Vert_{\dot H^{3/2}}\\
 \end{eqnarray*}
 
 \noindent {\bf{6.4. Estimate  of  $T_{1,4}$}}  \\

The control of $T_{1,4}$ is done thanks to the following decomposition which is analogous to the decomposition of $I_{1,3,2}$ (see \eqref{I132}), namely
\begin{eqnarray*}
T_{1,4}&=& - \frac{1}{2}\int \mathcal{H}\partial^{3}_{x}f \ \int \frac{f_{xx}(x)-f_{xx}(x-\alpha)}{\alpha} \  \ \int_{0}^{\infty} \delta  e^{-\delta}  
\frac{f_{x}(x+\alpha)+f_{x}(x-\alpha)-2f_{x}(x)}{\alpha} \\
&& \times \ \cos(\frac{\delta}{2}(\Delta_{\alpha}f-\bar\Delta_{\alpha}f)) \sin(\frac{\delta}{2}(\Delta_{\alpha}f+\bar\Delta_{\alpha}f))  \ d\delta \ d\alpha \ dx \\
&+& \frac{1}{2} \int \mathcal{H}\partial^{3}_{x}f \ \int \frac{f_{xx}(x)-f_{xx}(x-\alpha)}{\alpha} \  \ \int_{0}^{\infty} \delta  e^{-\delta}  
 \ \frac{\int_{0}^{\alpha} f_{x}(x-s)+f_{x}(x+s)-2f_{x}(x) \ ds}{\alpha^2} \\
 &&\times \cos(\frac{\delta}{2}(\Delta_{\alpha}f-\bar\Delta_{\alpha}f)) \sin(\frac{\delta}{2}(\Delta_{\alpha}f+\bar\Delta_{\alpha}f))  \ d\delta \ d\alpha \ dx \\
 &=& T_{1,4,1} + T_{1,4,2}.
\end{eqnarray*}

We have
\begin{eqnarray*}
\vert  T_{1,4,1} \vert &\lesssim& \Vert f \Vert_{\dot H^{3}} \int \frac{\Vert f_{xx}(x)-f_{xx}(x-\alpha) \Vert_{L^{\infty}}}{\vert\alpha \vert} \   
\frac{\Vert f_x(x+\alpha)+f_x(x-\alpha)-2f_x(x) \Vert_{L^{2}}}{\vert\alpha \vert} \ d\alpha  \\
&\leq&  \Vert f \Vert_{\dot H^{3}} \Vert f \Vert_{\dot B^{5/2}_{\infty,2}} \Vert f \Vert_{\dot B^{3/2}_{2,2}} \\
&\leq&  \Vert f \Vert^{2}_{\dot H^{3}} \Vert f\Vert_{\dot H^{3/2}}
\end{eqnarray*}
For   $T_{1,4,2}$, we write
\begin{eqnarray*}
\vert  T_{1,4,2} \vert &\leq&\frac{1}{2}\Vert f \Vert_{\dot H^{3}} \ \int \frac{\Vert f_{xx}(x)-f_{xx}(x-\alpha) \Vert_{L^{\infty}}}{\vert \alpha \vert^{3}} \  \ \int_{0}^{\infty} \delta  e^{-\delta} \\
&\times&   \vert \alpha \vert^{q+\frac{1}{\bar r}} \left(\int_{0}^{\alpha} \frac{ \Vert f_{x}(x-s)+f_{x}(x+s)-2f_{x}(x) \Vert^{r}_{L^{2}}}{\vert s\vert^{qr}} \ ds\right)^{1/r}  \frac{\Vert \delta_{\alpha}f+\bar\delta_{\alpha}f \Vert_{L^{\infty}}}{\vert\alpha\vert}  \ d\delta \ d\alpha  \\
 &\leq&\frac{\Gamma(3)}{2}\Vert f \Vert_{\dot H^{3}} \Vert f_{x} \Vert_{ \dot B^{q-\frac{1}{r}}_{2,r}} \ \int \frac{\Vert f_{xx}(x)-f_{xx}(x-\alpha) \Vert_{L^{\infty}}}{\vert \alpha \vert^{2-q-\frac{1}{\bar r}}} \  \frac{\Vert \delta_{\alpha}f+\bar\delta_{\alpha}f \Vert_{L^{\infty}}}{\alpha^2} \ d\alpha  \\
 &\lesssim& \Vert f \Vert_{\dot H^{3}} \Vert f_x\Vert_{ \dot B^{q-\frac{1}{r}}_{2,r}} \left(\int \frac{\Vert f_{xx}(x)-f_{xx}(x-\alpha) \Vert^{2}_{L^{\infty}}}{\vert \alpha \vert^{4-2q-\frac{2}{\bar r}}} \ d\alpha \int \frac{\Vert \delta_{\alpha}f+\bar\delta_{\alpha}f \Vert^{2}_{L^{\infty}}}{\alpha^{4}} \ d\alpha \right)^{1/2} \\
 &\lesssim& \Vert f \Vert_{\dot H^{3}} \Vert f \Vert_{ \dot B^{q+1-\frac{1}{r}}_{2,r}} \Vert f \Vert_{\dot B^{7/2-q-\frac{1}{\bar r}}_{\infty,2}}\Vert f \Vert_{\dot B^{3/2}_{\infty,2}} \\
 &\lesssim&\Vert f \Vert_{\dot H^{3}} \Vert f \Vert_{\dot H^{2}} \Vert f \Vert_{ \dot B^{q+\frac{1}{2}}_{2,2}} \Vert f \Vert_{\dot B^{3-q}_{\infty,2}}
 \end{eqnarray*}
where we have chosen $\bar r= r=2$ and then by choosing $q=1$  we obtain
\begin{eqnarray*}
\vert  T_{1,4,2} \vert & \lesssim& \Vert f \Vert_{\dot H^{3}} \Vert f \Vert_{\dot H^{3/2}} \Vert f \Vert_{\dot H^{5/2}}\Vert f \Vert_{\dot H^{2}} \\
 & \lesssim&  \Vert f \Vert^{2}_{\dot H^{3}} \Vert f \Vert^{2}_{\dot H^{3/2}}
\end{eqnarray*} \label{eq:croisés}
Therefore,
\begin{equation*}
\vert T_{1,4} \vert \lesssim \Vert f \Vert^{2}_{\dot H^{3}}(\Vert f \Vert_{\dot H^{3/2}}+\Vert f \Vert^{2}_{\dot H^{3/2}}).
\end{equation*}

\noindent {\bf{6.5. Estimate  of  $T_{1,5}$}}  \\

We now estimate  $T_{1,5}$, analogously to $I_{133}$ (see \eqref{I132}) we use the identity \eqref{formuleSa}. So that,
\begin{eqnarray*}
  T_{1,5} &=& -\frac{1}{2}\int \mathcal{H}\partial^{3}_{x}f  \int \frac{f_{xx}(x)-f_{xx}(x-\alpha)}{\alpha} \  \ \int_{0}^{\infty} \delta  e^{-\delta}  
  \sin(\frac{\delta}{2}(\Delta_{\alpha}f-\bar\Delta_{\alpha}f)) \partial_{\alpha} S  \cos(\frac{\delta}{2}(\Delta_{\alpha}f+\bar\Delta_{\alpha}f))  \ d\delta \ d\alpha \ dx \\
 &=& \frac{1}{2} \int \mathcal{H}\partial^{3}_{x}f  \int \frac{f_{xx}(x)-f_{xx}(x-\alpha)}{\alpha}  \int_{0}^{\infty} \delta  e^{-\delta}  
 \ \sin(\frac{\delta}{2}(\Delta_{\alpha}f-\bar\Delta_{\alpha}f)) \bar \Delta_{\alpha} f_{x} \ \cos(\frac{\delta}{2}(\Delta_{\alpha}f+\bar\Delta_{\alpha}f))  \ d\delta \ d\alpha \ dx \\
&+&\frac{1}{2} \int \mathcal{H}\partial^{3}_{x}f \ \int \frac{f_{xx}(x)-f_{xx}(x-\alpha)}{\alpha}    \int_{0}^{\infty} \delta  e^{-\delta} 
  \sin(\frac{\delta}{2}(\Delta_{\alpha}f-\bar\Delta_{\alpha}f)) \Delta_{\alpha}f_{x} \ \cos(\frac{\delta}{2}(\Delta_{\alpha}f+\bar\Delta_{\alpha}f))  \ d\delta \ d\alpha \ dx \\
 &-&\frac{1}{2} \int \mathcal{H}\partial^{3}_{x}f \ \int \frac{f_{xx}(x)-f_{xx}(x-\alpha)}{\alpha}  \int_{0}^{\infty} \delta  e^{-\delta}  
  \sin(\frac{\delta}{2}(\Delta_{\alpha}f-\bar\Delta_{\alpha}f)) \frac{f(x+\alpha)+f(x-\alpha)-2f(x)}{\alpha^2} \\
  && \times\cos(\frac{\delta}{2}(\Delta_{\alpha}f+\bar\Delta_{\alpha}f))  \ d\delta \ d\alpha \ dx \\
  &=& \sum_{j=1}^{3}    T_{1,5, j}.
  \end{eqnarray*}
 
 The estimate of $T_{1,5,1}$ is relatively easy, indeed, it suffices to write
  \begin{eqnarray*}
  \vert T_{1,5,1} \vert &\leq&\frac{\Gamma(2)}{2} \Vert f \Vert_{\dot H^{3}} \left(\int \frac{\Vert f_{xx}(x)-f_{xx}(x-\alpha) \Vert^{2}_{L^{\infty}}}{\vert \alpha \vert^{2}} \ d\alpha \int \frac{\Vert f_{x}(x)-f_{x}(x+\alpha) \Vert^{2}_{L^{2}}}{\vert \alpha \vert^{2}} \ d\alpha \right)^{1/2} \\
   & \lesssim& \Vert f \Vert_{\dot H^{3}} \Vert f \Vert_{\dot B^{5/2}_{\infty,2}} \Vert f \Vert_{\dot B^{3/2}_{2,2}} \\
       & \lesssim&  \Vert f \Vert^{2}_{\dot H^{3}} \Vert f \Vert_{\dot H^{3/2}}
  \end{eqnarray*}
  As well, one may easily estimate $T_{1,5,2}$ by writting
  \begin{eqnarray*}
   \vert T_{1,5,2} \vert  &\leq&\frac{\Gamma(2)}{2} \Vert f \Vert_{\dot H^{3}} \left(\int \frac{\Vert f_{xx}(x)-f_{xx}(x-\alpha) \Vert^{2}_{L^{\infty}}}{\vert \alpha \vert^{2}} \ d\alpha \int \frac{\Vert f_{x}(x)-f_{x}(x+\alpha) \Vert^{2}_{L^{2}}}{\vert \alpha \vert^{2}} \ d\alpha \right)^{1/2} \\
   & \lesssim&\Vert f \Vert_{\dot H^{3}} \Vert f \Vert_{\dot B^{5/2}_{\infty,2}} \Vert f \Vert_{\dot B^{3/2}_{2,2}} \\
       &\lesssim  &  \Vert f \Vert^{2}_{\dot H^{3}} \Vert f \Vert_{\dot H^{3/2}}
         \end{eqnarray*}
         For $T_{1,5,3}$, it suffices to write
  \begin{eqnarray*}
   \vert T_{1,5,3} \vert &\leq&\frac{\Gamma(2)}{2} \Vert f \Vert_{\dot H^{3}} \left(\int \frac{\Vert f_{xx}(x)-f_{xx}(x-\alpha) \Vert^{2}_{L^{\infty}}}{\vert \alpha \vert^{2}} \ d\alpha \int \frac{\Vert f(x+\alpha)+f(x-\alpha)-2f(x)\Vert^{2}_{L^{2}}}{\vert \alpha \vert^{4}} \ d\alpha \right)^{1/2} \\
   & \lesssim&\Vert f \Vert_{\dot H^{3}} \left(\int \frac{\Vert f_{xx}(x)-f_{xx}(x-\alpha) \Vert^{2}_{L^{\infty}}}{\vert \alpha \vert^{2}} \ d\alpha \int \frac{\Vert f(x+\alpha)+f(x-\alpha)-2f(x)\Vert^{2}_{L^{2}}}{\vert \alpha \vert^{4}} \ d\alpha \right)^{1/2} \\
   & \lesssim&\Vert f \Vert_{\dot H^{3}} \Vert f_{xx} \Vert_{\dot B^{1/2}_{\infty,2}}\Vert f \Vert_{\dot B^{3/2}_{2,2}} \\
   & \lesssim&\Vert f \Vert^{2}_{\dot H^{3}} \Vert f \Vert_{\dot H^{3/2}}
   \end{eqnarray*}
   Therefore,
   \begin{equation} \label{eq:croisés2}
    \vert T_{1,5} \vert \lesssim \Vert f \Vert^{2}_{\dot H^{3}} \Vert f \Vert_{\dot H^{3/2}}
   \end{equation}
  It remains to estimate $T_{1,6}$, this is the purpose of the next subsection. \\

\noindent {\bf{6.5. Estimate  of  $T_{1,6}$}}  \\

As we did for $I_{1,3,4}$ (see equality \eqref{I134}), we first rewrite $T_{1,6}$ in term of controlled commutators. 
\begin{eqnarray*}
T_{1,6} &=& - \int \mathcal{H}\partial^{3}_{x}f_{xx} \ \int \frac{\partial^{3}_{x}f(x)}{\alpha} \  \ \int_{0}^{\infty}   e^{-\delta}  
\sin(\frac{\delta}{2}(\Delta_{\alpha}f-\bar\Delta_{\alpha}f)) \sin(\frac{\delta}{2}(\Delta_{\alpha}f+\bar\Delta_{\alpha}f))  \ d\delta \ d\alpha \ dx \\
&=& \frac{1}{2} \int  \partial^{3}_{x}f \ \int_{0}^{\infty}     e^{-\delta} \int \frac{1}{\alpha}
 \left[\mathcal{H} , \sin(\frac{\delta}{2}(\Delta_{\alpha}f-\bar\Delta_{\alpha}f)) \sin(\frac{\delta}{2}(\Delta_{\alpha}f+\bar\Delta_{\alpha}f)) \right]\partial^{3}_{x}f   \ d\delta \ d\alpha \ dx \\
 &=&  \frac{1}{2}\int \int \frac{\partial^{3}_{x}f-\partial^{3}_{x}f(x-\alpha)}{\alpha} \ \int_{0}^{\infty}     e^{-\delta}  
 \left[\mathcal{H} , \sin(\frac{\delta}{2}(\Delta_{\alpha}f-\bar\Delta_{\alpha}f)) \sin(\frac{\delta}{2}(\Delta_{\alpha}f+\bar\Delta_{\alpha}f)) \right]\partial^{3}_{x}f   \ d\delta \ d\alpha \ dx \\
  &+&\frac{1}{2} \int\int \frac{\partial^{3}_{x}f(x-\alpha)}{\alpha} \ \int_{0}^{\infty}     e^{-\delta}  
 \left[\mathcal{H} , \sin(\frac{\delta}{2}(\Delta_{\alpha}f-\bar\Delta_{\alpha}f)) \sin(\frac{\delta}{2}(\Delta_{\alpha}f+\bar\Delta_{\alpha}f)) \right]\partial^{3}_{x}f   \ d\delta \ d\alpha \ dx.
   \end{eqnarray*}
   Finally, by integrating by parts one obtains
   \begin{eqnarray*}
T_{1,6,1} &=& -\frac{1}{2}\int\int \frac{f_{xx}(x)-f_{xx}(x-\alpha)}{\alpha} \ \int_{0}^{\infty}     e^{-\delta}  
 \partial_{x}\left[\mathcal{H} , \sin(\frac{\delta}{2}(\Delta_{\alpha}f-\bar\Delta_{\alpha}f)) \sin(\frac{\delta}{2}(\Delta_{\alpha}f+\bar\Delta_{\alpha}f)) \right]\partial^{3}_{x}f \  d\delta \ d\alpha \ dx \\
  &-&\frac{1}{2} \int \int \frac{f_{xx}(x-\alpha)-f_{xx}(x)}{\alpha^2} \ \int_{0}^{\infty}     e^{-\delta}  
 \left[\mathcal{H} , \sin(\frac{\delta}{2}(\Delta_{\alpha}f-\bar\Delta_{\alpha}f)) \sin(\frac{\delta}{2}(\Delta_{\alpha}f+\bar\Delta_{\alpha}f)) \right]\partial^{3}_{x}f   \ d\delta \ d\alpha \ dx \\
 &+&\frac{1}{2} \int \int \frac{f_{xx}(x-\alpha)-f_{xx}(x)}{\alpha} \ \int_{0}^{\infty}     e^{-\delta}  
 \partial_{\alpha}\left[\mathcal{H} , \sin(\frac{\delta}{2}(\Delta_{\alpha}f-\bar\Delta_{\alpha}f)) \sin(\frac{\delta}{2}(\Delta_{\alpha}f+\bar\Delta_{\alpha}f)) \right]\partial^{3}_{x}f   \ d\delta \ d\alpha \ dx \\
 &=& T_{1,6,1}+T_{1,6,2}+T_{1,6,3}.
   \end{eqnarray*}
 By using the generalized Calder\'on commutator estimate \eqref{cz}
 along with some classical Besov embeddings, we may  control  $T_{1,6,1}$ as follows. 
\begin{eqnarray*}
 \vert T_{1,6,1} \vert &\lesssim&  \Vert f\Vert_{\dot H^{3}} \int  \frac{\Vert f_{xx}(x)-f_{xx}(x-\alpha) \Vert_{L^{2}}}{\vert\alpha \vert} \frac{\Vert  f_{x}(x)-f_{x}(x-\alpha)\Vert_{L^{\infty}}}{\vert\alpha \vert}  \frac{\Vert f (x-\alpha)+f(x+\alpha)-2f(x) \Vert_{L^{\infty}}}{\vert\alpha \vert}  \ d\alpha \\
 &+& \Vert f\Vert_{\dot H^{3}} \int  \frac{\Vert f_{xx}(x)-f_{xx}(x-\alpha) \Vert_{L^{2}}}{\vert\alpha \vert} \frac{\Vert  f_{x}(x)-f_{x}(x+\alpha)\Vert_{L^{\infty}}}{\vert\alpha \vert}  \frac{\Vert f (x-\alpha)+f(x+\alpha)-2f(x) \Vert_{L^{\infty}}}{\vert\alpha \vert}  \ d\alpha \\
 &+&  \Vert f\Vert_{\dot H^{3}} \int  \frac{\Vert f_{xx}(x)-f_{xx}(x-\alpha) \Vert_{L^{2}}}{\vert\alpha \vert}   \frac{\Vert f_x (x-\alpha)+f_x (x+\alpha)-2f_x (x) \Vert_{L^{\infty}}}{\vert\alpha \vert}  \ d\alpha  \\
 &\lesssim&  \Vert f\Vert_{\dot H^{3}}   \left(\int \frac{\Vert f_{xx}(x)-f_{xx}(x-\alpha) \Vert^{2}_{L^{2}}}{\alpha^2} \ d\alpha \right)^{1/2} \left(\int \frac{\Vert f (x-\alpha)+f(x+\alpha)-2f(x) \Vert^{4}_{L^{\infty}}}{\vert\alpha \vert^{5}}  \ d\alpha \right)^{1/4}  \\
 && \times \left[\left(\int\frac{\Vert  f_{x}(x)-f_{x}(x-\alpha)\Vert^{4}_{L^{\infty}}}{\vert\alpha \vert^3} \ d\alpha \right)^{1/4}+ \left(\int\frac{\Vert  f_{x}(x)-f_{x}(x+\alpha)\Vert^{4}_{L^{\infty}}}{\vert\alpha \vert^3} \ d\alpha \right)^{1/4}\right]\\
 &+&  \Vert f\Vert_{\dot H^{3}}   \left(\int \frac{\Vert f_{xx}(x)-f_{xx}(x-\alpha) \Vert^{2}_{L^{2}}}{\alpha^2} \ d\alpha \right)^{1/2} \left(\int \frac{\Vert f_{x} (x-\alpha)+f_{x}(x+\alpha)-2f_{x}(x) \Vert^{2}_{L^{\infty}}}{\alpha^{2}}  \ d\alpha \right)^{1/2}  \\
 &\lesssim&  \Vert f\Vert_{\dot H^{3}} \Vert f \Vert_{\dot H^{5/2}}\left( \Vert f \Vert_{\dot B^{1}_{\infty,4}} \Vert f \Vert_{\dot B^{3/2}_{\infty,4}}  + \Vert f \Vert_{\dot B_{\infty,2}^{3/2}}\right) \\
 & \lesssim&  \Vert f\Vert_{\dot H^{3}} \Vert f \Vert_{\dot H^{5/2}}\left( \Vert f \Vert_{\dot B^{1}_{\infty,4}} \Vert f \Vert_{\dot B^{3/2}_{\infty,4}}  + \Vert f \Vert_{\dot B_{\infty,2}^{3/2}}\right) \\
 &\lesssim&  \Vert f\Vert_{\dot H^{3}} \Vert f \Vert^{1/3}_{\dot H^{3/2}} \Vert f \Vert^{2/3}_{\dot H^{3}}\left( \Vert f \Vert_{\dot H^{3/2}} \Vert f \Vert_{\dot B^{2}_{2,2}}  + \Vert f \Vert_{\dot H^{2}}\right) \\
 &\lesssim&  \Vert f\Vert^{2}_{\dot H^{3}}  \Vert f \Vert^2_{\dot H^{3/2}} + \Vert f\Vert^{2}_{\dot H^{3}}  \Vert f \Vert_{\dot H^{3/2}}. 
 \end{eqnarray*}
Then, we estimate $T_{1,6,2}$. Analogously to $I_{1,3,4,2}$ (see \eqref{I1342}) we see that this term may be rewritten as
\begin{eqnarray*}
T_{1,6,2}&=& \frac{1}{2}\int \mathcal{H}\partial^{3}_{x}f \int \frac{f_{xx}(x-\alpha)-f_{xx}(x)}{\alpha^2} \ \int_{0}^{\infty}     e^{-\delta}  \sin(\frac{\delta}{2}(\Delta_{\alpha}f-\bar\Delta_{\alpha}f)) \sin(\frac{\delta}{2}(\Delta_{\alpha}f+\bar\Delta_{\alpha}f))  \  \ d\delta \ d\alpha \ dx \\
&+& \frac{1}{2}\int \partial^{3}_{x}f  \int \frac{\mathcal{H}f_{xx}(x-\alpha)-\mathcal{H}f_{xx}(x)}{\alpha^2} \ \int_{0}^{\infty}     e^{-\delta}  \sin(\frac{\delta}{2}(\Delta_{\alpha}f-\bar\Delta_{\alpha}f)) \sin(\frac{\delta}{2}(\Delta_{\alpha}f+\bar\Delta_{\alpha}f))  \  \ d\delta \ d\alpha \ dx \\
&=& T_{1,6,2,1}+ T_{1,6,2,2}.
\end{eqnarray*}
To estimate $T_{1,6,2,1}$, we use the embedding $\dot H^2 \hookrightarrow \dot B^{3/2}_{\infty,2}$ along with classical interpolation inequalites to get
\begin{eqnarray*}
\vert T_{1,6,2,1} \vert &\lesssim& \Vert f \Vert_{\dot H^{3}} \int \frac{\Vert f_{xx}(x-\alpha)-f_{xx}(x) \Vert_{L^{2}}}{\alpha^2} \ \int_{0}^{\infty}   \delta  e^{-\delta}  \frac{\Vert \delta_{\alpha}f+\bar\delta_{\alpha}f \Vert_{L^{\infty}}}{\vert\alpha \vert}  \  \ d\delta \ d\alpha   \\
&\lesssim&  \Vert f \Vert_{\dot H^{3}} \left(\int \frac{\Vert f_{xx}(x-\alpha)-f_{xx}(x) \Vert^{2}_{L^{2}}}{\alpha^2} \ d\alpha \right)^{1/2}  \left( \int \frac{\Vert \delta_{\alpha}f+\bar\delta_{\alpha}f \Vert^{2}_{L^{\infty}}}{\alpha^4} \ d\alpha \right)^{1/2} \\
&\lesssim& \Vert f \Vert_{\dot H^{3}} \Vert f \Vert_{\dot H^{5/2}}  \Vert f \Vert_{\dot H^{2}} \\
&\lesssim&\Vert f \Vert^{2}_{\dot H^{3}}  \Vert f \Vert_{\dot H^{3/2}}
\end{eqnarray*}
Since  $\mathcal{H}$ is continuous on $L^{2}$, one gets the same control for $T_{1,6,2,2}$. Therefore,
\begin{eqnarray*}
\vert T_{1,6,2} \vert \lesssim  \Vert f \Vert^{2}_{\dot H^3} \Vert f \Vert_{\dot H^{3/2}}
\end{eqnarray*}
It remains to estimate $T_{1,6,3}$, to do so, we use the following  decomposition (which is analogous to equality \eqref{I1343}).
\begin{eqnarray*}
 T_{1,6,3}&=&-\frac{1}{4}\int \partial_{x}^{3} f \int  \frac{\mathcal{H}f_{xx}(x-\alpha)-\mathcal{H}f_{xx}(x)}{\alpha} \ \int_{0}^{\infty}   \delta  e^{-\delta}  \partial_{\alpha}(\Delta_{\alpha}f-\bar\Delta_{\alpha}f) \ \cos(\frac{\delta}{2}(\Delta_{\alpha}f-\bar\Delta_{\alpha}f)) \\
 && \times \sin(\frac{\delta}{2}(\Delta_{\alpha}f+\bar\Delta_{\alpha}f))  \  \  \ d\delta \ d\alpha \ dx \\
 &-&\frac{1}{4}\int \mathcal{H}\partial_{x}^{3}f  \int \frac{f_{xx}(x-\alpha)-f_{xx}(x)}{\alpha} \ \int_{0}^{\infty}   \delta  e^{-\delta}  \partial_{\alpha}(\Delta_{\alpha}f-\bar\Delta_{\alpha}f) \ \cos(\frac{\delta}{2}(\Delta_{\alpha}f-\bar\Delta_{\alpha}f)) \\
 && \times \sin(\frac{\delta}{2}(\Delta_{\alpha}f+\bar\Delta_{\alpha}f)  \  \  \ d\delta \ d\alpha \ dx \\
 &-&\frac{1}{4} \int \partial_{x}^{3}f \int \frac{\mathcal{H}f_{xx}(x-\alpha)-\mathcal{H}f_{xx}(x)}{\alpha} \ \int_{0}^{\infty}   \delta  e^{-\delta}  \partial_{\alpha}(\Delta_{\alpha}f+\bar\Delta_{\alpha}f) \ \cos(\frac{\delta}{2}(\Delta_{\alpha}f+\bar\Delta_{\alpha}f)) \\
 && \times \sin(\frac{\delta}{2}(\Delta_{\alpha}f-\bar\Delta_{\alpha}f)  \  \ d\delta \ d\alpha \ dx \\
 &-&\frac{1}{4} \int \partial_{x}^{3} f \int \frac{f_{xx}(x-\alpha)-f_{xx}(x)}{\alpha} \ \int_{0}^{\infty}   \delta  e^{-\delta}  \partial_{\alpha}(\Delta_{\alpha}f+\bar\Delta_{\alpha}f) \ \cos(\frac{\delta}{2}(\Delta_{\alpha}f+\bar\Delta_{\alpha}f)) \\
 && \times \sin(\frac{\delta}{2}(\Delta_{\alpha}f-\bar\Delta_{\alpha}f)  \  \ d\delta \ d\alpha \ dx. 
 \end{eqnarray*}
Since we shall do $L^{p}$ estimates, $p\in(1,\infty)$ on the terms involving $\mathcal{H}$. We observe that these terms has the same regularity as the terms $I_{1,3,4,3}$ (see \eqref{I1343}) and therefore one anagolously infers that,
$$
\left \vert  T_{1,6,3} \right\vert \lesssim \Vert f \Vert^{2}_{\dot H^{3}} \Vert f \Vert^{2}_{\dot H^{3/2}} + \Vert f \Vert^{2}_{\dot H^{3}} \Vert f \Vert_{\dot H^{3/2}},
$$
hence,
$$
\left \vert T_{1,6} \right\vert \lesssim \Vert f \Vert^{2}_{\dot H^{3}} (\Vert f \Vert_{\dot H^{3/2}}+\Vert f \Vert^{2}_{\dot H^{3/2}}).
$$

\noindent Therefore, we  get that
 
 $$\vert T_{1} \vert \lesssim \Vert f \Vert^{2}_{\dot H^{3}}   (\Vert f \Vert^{2}_{\dot H^{3/2}} + \Vert f \Vert_{\dot H^{3/2}})  - \pi \Vert f \Vert^{2}_{\dot H^{3}} + \pi\frac{K^2}{1+K^2}\Vert  f \Vert^{2}_{\dot H^{3}}.$$
 
\noindent  Finally, we obtain that,
 
 \begin{equation} \label{T1}
 \vert T_{1} \vert \lesssim \Vert f \Vert^{2}_{\dot H^{3}}   (\Vert f \Vert^{2}_{\dot H^{3/2}} + \Vert f \Vert_{\dot H^{3/2}})  - \frac{\pi}{1+K^2}\Vert  f \Vert^{2}_{\dot H^{3}}.
\end{equation}

\noindent {\bf{6.7. Estimate  of  $T_2$}}  \\

  Recall that,
  \begin{eqnarray*}
 T_{2}=  -3\int  \mathcal{H}\partial_{x}^{3}f \ \int \partial_{xx} \Delta_{\alpha} f  \ \partial_{x} \Delta_{\alpha} f  \ \int_{0}^{\infty} \delta  e^{-\delta}  \sin(\delta\Delta_{\alpha}f(x) ) \ d\delta \ d\alpha \ dx.       
  \end{eqnarray*}
 To estimate this term, it suffices to observe that,
  $$
  \vert T_{2} \vert \lesssim \Vert f\Vert_{\dot H^{3}} \left(\int \frac{\Vert f_{xx}(x)-f_{xx}(x-\alpha) \Vert^{2}_{L^{2}}}{\alpha^{2}} \ d\alpha\right)^{1/2} \left(\int \frac{\Vert f_{x}(x)-f_{x}(x-\alpha) \Vert^{2}_{L^{\infty}}}{\alpha^{2}} \ d\alpha\right)^{1/2}.
  $$
Therefore,
  $$
  \vert T_{2} \vert \lesssim \Vert f \Vert_{\dot H^{3}} \Vert f \Vert_{\dot H^{5/2}} \Vert f_{x}\Vert_{\dot B^{1/2}_{\infty,2}}.
  $$
  Hence, using the embedding $\dot H^{1} \hookrightarrow \dot B^{1/2}_{\infty,2}$ along with Sobolev interpolations $\dot H^{2}=\left[\dot H^{3/2}, \dot H^{3}\right]_{\frac{2}{3},\frac{1}{3}}$ and $\dot H^{5/2}=\left[\dot H^{3/2}, \dot H^{3}\right]_{\frac{1}{3},\frac{2}{3}}$, one finds 
  \begin{equation} \label{T2}
  \vert T_{2} \vert \lesssim \Vert f \Vert^{2}_{\dot H^{3}} \Vert f \Vert_{\dot H^{3/2}}
  \end{equation}
\noindent {\bf{6.8. Estimate  of  $T_3$}}  \\

It suffices to observe for instance that $ \dot B^{5/3}_{6,3}\hookleftarrow\dot H^{2}=\left[\dot H^{3/2}, \dot H^{3}\right]_{\frac{2}{3},\frac{1}{3}}$, so that
 \begin{eqnarray} \label{T3}
T_3  &=&\nonumber- 2\int \mathcal{H}\partial^{3}_{x} f 
\ \int (\partial_{x} \Delta_{\alpha}f)^3 \int_{0}^{\infty} \delta^{2} e^{-\delta}  \ \cos(\delta \Delta_{\alpha} f)  \ dx \ d\alpha \ d\delta \\
&\lesssim&\nonumber \Vert f \Vert_{\dot H^{3}} \Vert f \Vert^{3}_{\dot B^{5/3}_{6,3}} \\
&\lesssim&\nonumber \Vert f \Vert^2_{\dot H^{3}} \Vert f \Vert^{2}_{\dot H^{3/2}} \\
\end{eqnarray}

\noindent Finally, by \eqref{T1}, \eqref{T2} and \eqref{T3} and by integrating in time $s\in[0,T]$ we find that for all $T>0$

\begin{equation*} 
 \Vert f \Vert^{2}_{\dot H^{5/2}}(T) + \frac{\pi}{1+K^2} \int_{0}^T \Vert f \Vert^{2}_{\dot H^{3}} \ ds \lesssim \Vert f_{0} \Vert_{\dot H^{5/2}} +  P(f) \int_{0}^T \Vert f \Vert^{2}_{\dot H^{3}} \ ds
\end{equation*}
where  $K=\displaystyle\sup_{(x,t)\in \mathbb R \times  [0,T]}\vert f_x(x,t)\vert$ and $P(f)=\Vert f \Vert_{L^{\infty}([0,T],\dot H^{3/2})} + \Vert f \Vert^{2}_{L^{\infty}([0,T],\dot H^{3/2})}$.  This ends the proof of Lemma \ref{ineg2} is proved. \\

\qed

\section{Proofs of Theorems \ref{res} and \ref{reg}}

We consider the following approximated system,
\begin{equation} 
\ (\tilde{\mathcal{M}_{\epsilon}}) \ : \\\left\{
\aligned
& \partial_{t}f_{\epsilon} (t,x) - \frac{\rho}{\pi} \int \partial_{x}\Delta_{\alpha}f_{\epsilon} \ \int_{0}^{\infty} e^{-\delta} \cos(\delta \Delta_{\alpha} f_{\epsilon}) \ d\delta \ d\alpha -\epsilon \Delta f_{\epsilon} =0\\
& f_{0,\epsilon}(x)=f_{0}(x)*\phi_{\epsilon}(x),
\endaligned
\right.
\end{equation}
where $\phi_\eps$ is a classical mollifier that is $\phi_{\epsilon}(x)=\eps^{-1} \phi(x\eps^{-1})$, $\phi \in \mathcal{D}(\mathbb R)$   so that $\phi$ is nonnegative and $\int\phi(x) dx =1$. In order to simplify the proof we assume that the initial data is of finite energy. If the Lipschitz norm remains bounded on a time interval $[0,T_{\eps})$ and if the $\Vert f_{0,\eps} \Vert_{H^{3/2}}$ is smaller than a constant that depends only on the $L^{\infty}([0,T_{\eps}), \dot W^{1,\infty})$ norm then we may locally solve the equation. Using the same {\it{a priori}} estimate as proved in section \ref{h32} for the regularized equation $(\tilde{\mathcal{M}_{\epsilon}})$ we may show that actually $T_{\eps}=+\infty$. One has  (uniformly in $\epsilon$)  that
$$
\Vert f_{0,\epsilon} \Vert_{H^{3/2}} \lesssim \Vert f_{0} \Vert_{H^{3/2}}
$$
as well as,
$$
\Vert f_{0,\epsilon} \Vert_{\dot W^{1,\infty}} \lesssim \Vert f_{0} \Vert_{\dot W^{1,\infty}}.
$$
Therefore, the regularized initial data converges strongly in $H^{3/2} \cap \dot W^{1,\infty}$. Let $\phi \in \mathcal{D}([0,T]\times\mathbb R)$ be nonnegative, from the {\it{a priori}} estimates we know that $\phi f_{\eps}$ is bounded in $L^{\infty}([0,T]; H^{3/2}\cap \dot W^{1,\infty}) \cap L^{2}([0,T]; H^{2})$. Since those spaces a separable Banach spaces,  thanks to the Banach-Alaoglu theorem, we may extract  from these solutions $f_{\eps}$ a subsequence $\{f_{\eps_{k}}\}_{k\geq0}$  that converges weakly to a solution $f \in  L^{2}([0,T]; H^{2})$ and *-weakly in $L^{\infty}([0,T]; H^{3/2})$. In order to obtain the convergence in the sense of distribution, one needs to prove a strong convergence in $(L^{2}L^2)_{loc}$. To do so, we need for instance to get a nice bound  on $\partial_{t} f_{\eps}$ locally in space and time, namely on 
$$
-\eps \Lambda^{2} f_{\eps}-\Lambda f_{\eps} -2\int \partial_{x} \Delta_{\alpha} f_{\eps} \int_{0}^{\infty} \delta e^{-\delta} \sin^{2}(\frac{\delta}{2}\Delta_{\alpha} f_{\eps})  \ d\delta \ d\alpha
$$

Since  $f_{\eps}$ is bounded in $L^{2}([0,T]; \dot H^{2})$ then  $\partial_{x}f_{\eps}$ is bounded in $L^{2}([0,T]; \dot H^{1})$. It is not difficult to see that the contribution coming from the nonlinear part of $\partial_{t} f_{\eps}$ is a locally bounded sequence in  $L^{2}H^{-1/4}$. Indeed, by using the dual form of the Sobolev embedding it suffices to have a bound on the $L^{2}L^{4/3}$ norm of the nonlinearity. By controling the product, we find that if $f_{\eps} \in \dot B^{3/2}_{2,2} \cap \dot B^{3/4}_{8,2}$ then the nonlinear part of $\partial_{t}f$ is bounded in $L^{2}\dot H^{-1/4}$. We know that $f_{\eps}$ is locally bounded in $L^{2}\dot H^{3/2} \cap L^2\dot H^{9/8}$. Since this space embedds into $\dot B^{3/2}_{2,2} \cap \dot B^{3/4}_{8,2}$, therefore, we get that the nonlinear term of $\partial_{t}f_{\epsilon}$ is a bounded sequence in $L^{2}\dot H^{-1/4}$. Since the linear part is locally  bounded in $L^{2}\dot H^1$,  we may use the Rellich compactness theorem (see e.g. \cite{PGLR}) to get the strong convergence of a subsequence in $(L^2L^2)_{loc}$. Consequently, the nonlinear term converges in $\mathcal{D}Õ$. It is then classical to prove that the limit is a solution of the equation. For higher regularity data (that is $f_{0} \in H^{5/2}$, with $f_{0}$ small enough in $\dot H^{3/2}$), using the {\emph{a priori}} estimates proved in section \ref{h52} we know that up to an extraction, $f_{\eps}$ is a bounded sequence in $L^{\infty}([0,T], H^{5/2})$ (in particular the  semi-norm $\sup_{t>0} \Vert f_{x} \Vert_{L^{\infty}}$ will remain bounded). Hence, we can pass to the weak limit as well since Rellich  gives also the strong compactness   in $(L^2L^2)_{loc}$. This concludes the result. \\

By doing the same {\it{a priori}} estimates for the difference of two solutions in $H^{3/2}$ (resp $H^{5/2}$), one observes that the uniqueness follows easily from the regularizing effect together with the fact that the $L^{\infty}_t H^{3/2}$ norm decays (resp $L^{\infty}_t H^{5/2}$). Hence, by Gr\"onwall's inequality gives the uniqueness in the usual way and we therefore omit the details.

\qed

\section*{Acknowledgments}
%%%%%%%%%%%%%%%
\noindent  Both D. C and O. L  were supported by the National Grant MTM2014-59488-P from the Spanish government and  ICMAT Severo Ochoa project SEV-2015-55. D.C. was partially supported by the ERC
Advanced Grant 788250.  O. L was supported by the Marie-Curie Grant, acronym: TRANSIC, from the FP7-IEF program, and the ERC through the Starting Grant project H2020-EU.1.1.-63922.

 \end{document}